\documentclass[smallextended]{svjour3}      

\usepackage{fullpage}
\usepackage{times}
\usepackage{ifthen}   
\usepackage{graphicx,graphics}
\usepackage{sidecap}
\usepackage[]{rotating}
\usepackage{setspace}
\usepackage{amsfonts}
\usepackage{color}
\usepackage{multirow}
\usepackage{longtable} 	
\usepackage{pbox}
\usepackage{datetime}
\usepackage{xspace}
\usepackage{enumitem} 
\usepackage[it,small]{caption}
\usepackage[T1]{fontenc}
\usepackage{doi}
\usepackage[numbers,sort]{natbib}
\usepackage{algorithmic}
\usepackage[linesnumbered,lined,boxed]{algorithm2e} 
\usepackage[misc]{ifsym}

\usepackage{academicons}
\def\orcid#1{\kern .1em\href{https://orcid.org/#1}{\includegraphics[keepaspectratio,width=1em]{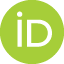}}}

\usepackage{booktabs}

\graphicspath{{./images/}} 
\usepackage{pgfgantt}

\usepackage{titlesec}
  \usepackage{etoolbox}
  \makeatletter
  \patchcmd{\ttlh@hang}{\parindent\z@}{\parindent\z@\leavevmode}{}{}
  \patchcmd{\ttlh@hang}{\noindent}{}{}{}
  \makeatother

\usepackage{mdframed}                 
\usepackage{colortbl}
\usepackage[most]{tcolorbox}
\usepackage{enumitem}

 \newcommand{\R}{\mathbb R}           
\newcommand{\dps}{\displaystyle}
\newcommand{\maxi}{\mathop{\mbox{maximize}}}
\newcommand{\mini}{\mathop{\mbox{minimize}}}

\newcommand{\st}{\mbox{subject to}}

\newcommand{\bvec}{\left(\begin{array}{c}}
\newcommand{\evec}{\end{array}\right)}

\newtheorem{assumption}{Assumption}

\usepackage{fancyhdr}
\usepackage[normalem]{ulem} 
 
{}

\usepackage{adjustbox}

\usepackage{xpatch}
\xpretocmd{\eqref}{Eq.~}{}{}

\begin{document}

\title{Modeling Design and Control Problems Involving Neural Network Surrogates}
\titlerunning{Optimization over Neural Network Surrogates}        

\author{
  Dominic Yang
  \and Prasanna Balaprakash
  \and Sven Leyffer
}

\institute{
  D. Yang \orcid{0000-0002-9453-2299} \Letter  \at
  University of California at Los Angeles, Los Angeles, USA\\
  \email{domyang@math.ucla.edu}
  \and
  P. Balaprakash \orcid{0000-0002-0292-5715} \at
  Argonne National Laboratory, Lemont, USA\\
  \email{pbalabra@anl.gov}
  \and
  S. Leyffer \orcid{0000-0001-8839-5876}
 \at
  Argonne National Laboratory, Lemont, USA\\
  \email{leyffer@anl.gov}
}

\date{Received: \today / Accepted: date}

\maketitle

\begin{abstract}%
  We consider nonlinear optimization problems that involve surrogate models represented by neural networks.
  We demonstrate first how to directly embed neural network evaluation into optimization models,
  highlight a difficulty with this approach that can prevent convergence, and then characterize stationarity of such models.
  We then present two alternative formulations of these problems in the specific case of feedforward neural networks with ReLU activation: as a mixed-integer optimization problem and as
  a mathematical program with complementarity constraints. 
  For the latter formulation we prove that stationarity at a point for this problem corresponds to stationarity of the embedded formulation.
  Each of these formulations may be solved with 
  state-of-the-art optimization methods, and we show how to obtain good initial feasible solutions for these methods.
  We compare our formulations on three practical applications arising in the design and control of 
  combustion engines, in the generation of adversarial attacks on classifier networks, and in the determination of optimal flows in an oil well network.%
  \subclass{MSC codes}
\end{abstract}

\begin{acknowledgements}
         {This material is based upon work supported by the 
U.S. Department of Energy,
Office of Science, 
Office of Advanced Scientific Computing Research, 
under Contract DE-AC02-06CH11357. This work was also supported by
the U.S. Department of Energy through grant DE-FG02-05ER25694. The first author was also supported through an NSF-MSGI fellowship.}
\end{acknowledgements}

\section{Introduction and Background}\label{sec:intro}

We are interested in solving general optimization problems that include deep neural networks (DNNs) that are used as surrogate models of complex functions (e.g., physical processes \cite{aithal2019maltese}, classification schemes \cite{krizhevsky2012imagenet, simonyan2014very, he2016deep}). 
In particular, we consider optimization problems
of the form
\begin{equation}\label{eq:OptiML}
    \mini_x \; f(\text{DNN}(x),x) \quad \st \; c(\text{DNN}(x),x) \le 0, \; x \in X,
\end{equation}
where $f$ and $c$ are smooth functions representing the objective function and constraints, respectively; $\text{DNN}(x)$ is the output of a DNN at $x$ (which we assume 
to have been previously trained on suitable data); $x$ are the optimization variables; and
$X \subset \R^n$ is a compact set that may include integer restrictions. We consider feedforward neural networks in this paper that are composed of a sequence of multiple layers of neurons. The values of neurons in layer $\ell$, $x^{\ell}$ are a nonlinear function (the \textit{activation function}) applied to a linear transformation of the values in the prior layer:
\begin{equation}\label{eq:DNNDef}
    x^{\ell} = \sigma(W^{{\ell}^T}x^{\ell-1} + b^{\ell}), \; \ell=1,\ldots,L.
\end{equation}
We have $x^0 = x$, $x^L = \text{DNN}(x)$, 
$L$ is the number of layers, $W$ and $b$ are weights determined by a training procedure, and $\sigma$ is the activation function applied componentwise. Unless otherwise specified, we will take $\sigma$ to be the ReLU function: $\text{ReLU(x)} = \max\{x, 0\}$.

An example of an optimization problem we  wish to solve is to minimize the output of a neural network regressor that predicts the quantity of emissions from automobile engine specifications.
In this case $X$ would be the set of existing automobile specifications, the constraints $c$ would be constraints to ensure the engine is realistic, and the objective would be a function of the emissions or engine performance. 
We will discuss this problem in  detail in Section \ref{sec:engine}.

To arrive at a tractable form of \eqref{eq:OptiML}, we make the following assumptions.
\begin{assumption}\label{ass:general}
We assume that the following conditions are satisfied for problem \eqref{eq:OptiML}:
\begin{enumerate}
    \item The objective function, $f: \R^n \to \R$, and the constraint functions, $c :\R^n \to \R^m$,
    are twice continuously differentiable {\em convex} functions.
    \item The set $X \subset \R^n$ is convex and compact.
\end{enumerate}
\end{assumption}
The most restrictive assumption is the convexity assumption. We can relax this assumption by leveraging 
standard global optimization techniques (see, e.g., \cite{sahinidis:96,tawarmalani.sahinidis:02,belotti2009couenne}),
at the expense of making the reformulated problem harder to solve. The smoothness assumption on $f$ and $c$
can be relaxed to Lipschitz continuity by using subgradients, and the compactness assumption is typically
satisfied as long as $x$ are constrained by bounds or the constraint function has compact level sets.
Below, we consider deep neural networks with ReLU activation functions that satisfy this assumption.

\paragraph{Outline and Contributions.}
We start by discussing three applications that employ DNN surrogates within an
optimization problem such as \eqref{eq:OptiML}:  (a) the minimization of emissions
in an engine design problem (as in \cite{aithal2019maltese}), (b) the generation of optimal adversarial examples to
``fool'' a given classifier (as in \cite{fischetti2018deep}), and (c) the optimization of the pump configuration for oil (as in \cite{grimstad2019relu}).
We develop a warmstart heuristic for both approaches that generates good initial guesses from the training data and helps us  overcome the challenges of the nonconvex formulation.
We also show how to add constraints that restrict the optimization problem to search only near where
existing training points can be added. 

Next, we demonstrate empirically that simply including a ReLU DNN within an optimization
problem can lead to poor convergence results. In particular, we have developed a new
nonlinear constraint for JuMP \cite{JuMP} that allows us to directly include DNNs within
an optimization model specified in JuMP. 
We then give a compact characterization of stationarity of the embedded formulation of our optimization problem.

Next, we consider alternative formulations and show that DNNs that use purely ReLU activation functions can be formulated as mixed-integer sets, building on \cite{fischetti2018deep,anderson2019strong}. We can then formulate \eqref{eq:OptiML} {\em equivalently} 
 as a convex mixed-integer nonlinear program (MINLP), which we refer to as the mixed-integer program (MIP) formulation. We also introduce a new formulation that results in a nonconvex nonlinear program (NLP) with complementarity constraints. We refer to this as the mathematical program with complementary constraints (MPCC) formulation. 
We prove this formulation has stationarity conditions equivalent to the embedded formulation. 
 Our reformulations involve a lifting into a higher-dimensional space in which 
the MINLP problem is convex. 

We show empirically for each of our applications that for moderately sized machine learning (ML) models
the resulting programs can be solved by using state-of-the-art MINLP and NLP solvers, making optimization problems with ML models tractable in practice.  We demonstrate that using the MIP formulation, we can find provably optimal solutions for small problems. Using the embedded and MPCC formulation, we show that we can address significantly larger networks at the cost of guaranteeing only locally optimal solutions, and we showcase scenarios where the MPCC formulation outperforms the embedded formulation  in terms of both optimal value found and consistency in convergence to a solution.

Throughout this paper we assume that the deep neural network has been trained
and is fixed for the optimization, and we do not consider the question of updating the neural network weights during the
optimization loop. One limitation of our approach is that we use standard MINLP solvers to 
tackle the reformulated MINLP, which limits the size of the neural network, $\text{DNN}(x)$, that can be used in the 
optimization. 

\paragraph{Related Work.}
Prior approaches to optimization over neural networks using MIP formulations include \citep{fischetti2018deep, cheng2017maximum, dutta2018output, khalil2018combinatorial, serra2020empirical, tjeng2017evaluating}. 
These approaches primarily model the piecewise ReLU constraint using standard big-M modeling tricks. 
They generally use the same basic formulation, but each augments the solve  by adding methods to tighten the big-M constraints \citep{fischetti2018deep, tjeng2017evaluating, grimstad2019relu}, decomposing the problem into smaller problems \citep{khalil2018combinatorial},  or adding local search routines \citep{dutta2018output}. 
Anderson et al.~\citep{anderson2019strong} provide an in-depth overview of how to strengthen these models to an ideal formulation with exponentially many constraints, as well as a method to separate in linear time.

Some other approaches to these problems using methods from MINLP have been tried. Cheon~\citep{cheon2020outer} solves an inverse
problem over ReLU constraints using an outer-approximation-inspired method, but without global optimality guarantees. 
Katz et al.~\citep{katz2017reluplex} use an approach from satisfiability modulo theory to address an optimization problem over neural networks. Scheweidtmann and Mitsos~\citep{schweidtmann2019deterministic} use a MINLP approach involving McCormick relaxations. 

One major focus of these optimization problems is on testing the resilience of neural networks against adversarial attack \citep{carlini2016evaluating}. This involves either  maximizing a notion of resilience \citep{cheng2017maximum} or finding minimal perturbations needed to misclassify an image \citep{fischetti2018deep}. Some work has also been done in using optimization to visualize features corresponding to neurons \citep{fischetti2018deep} and for  surrogate optimization  
in the context of optimizing the production of a set of oil wells \citep{grimstad2019relu}.
Other applications involve the use of neural networks as surrogates in the context of policy design for reinforcement learning \cite{ryu2019caql, delarue2020reinforcement}. A recent paper by Papalexopoulos et al.~\citep{papalexopoulos2021constrained} discusses the use of ReLU neural network surrogates for the purposes of black-box optimization.


\section{Modeling Optimization Applications Involving Neural Network Surrogates}\label{sec:model}

In this section we describe three optimization models that make use of neural-network
surrogates, and we discuss some of the challenges that arise.

\subsection{Optimal Design of Combustion Engine}\label{sec:engine}

Automobile engine operation is typically modeled by using highly intensive physics-based simulation code, which even on powerful computers can take hours just to model even a short drive. Hence, modeling the evaluation of this simulation code by using a neural network surrogate model can produce significant time savings, at the cost of producing slightly less accurate results. A study examining this approach is documented in \citep{aithal2019maltese}.

Given a trained neural network, we could formulate several optimization problems answering questions related to the operation of an engine on a given commute. The problem we consider is that  of minimizing emissions over the course of a given drive. The resulting solution will then be both the engine type and the drive style (e.g., RPM at all times in the drive) that produces the most environmentally efficient commute. The following section demonstrates how such a problem can be formulated.

\subsubsection{Simple Engine Design and Control Problem}
Suppose we have a trained DNN that predicts engine behavior based on engine specifications and driving parameters, given in Table \ref{tab:engineDNN}. The DNN is trained by using 64 trips each split into 1,500 observations (25 minutes observed at second intervals).

\begin{table}[htbp]
    \centering
    \begin{tabular}{l||l}
    Input Parameter & Output Parameter \\ \hline
    fuel injection (g) /s      &  nitrogen oxide, NO /s \\
    engine RPM     /s      &  carbon monoxide, CO /s  \\
    compression ratio     &  torque /s \\
    \end{tabular}
\caption{Input/output parameters of engine DNN model. Time-dependent parameters are shown with time per second (/s).}
\label{tab:engineDNN}
\end{table}

The problem we propose to solve is the optimal design and control of an engine for a given 25-minute trip. We will use a prescribed torque profile as a surrogate for the trip characteristics and minimize a weighted sum of NO and CO output. We have the following design variables  that are input to the DNN: fuel injection, $f$; engine RPM, $r$; and compression ratio, $c$. Note that $f$ and $r$ are characteristics that change over time whereas $c$ is an engine parameter that is fixed for the full drive. We also have variables NO and CO that represent nitrogen oxide and carbon monoxide emissions produced by the engine, as well as a variable torque that indicates engine torque. Each of these quantities is predicted by the neural network for each time interval of the drive. Formally, we state  the following optimal design and control problem over the time horizon $t=1,\ldots,T$ (where $T$ is the number of time intervals):
\begin{equation}\label{eq:optEngine}
    \begin{array}{lll}
        \dps \mini_{f,r,c,\text{CO},\text{NO},\text{Torque}} & \dps \sum_{t=1}^T (\text{NO}_t + \lambda \text{CO}_t)\Delta t \qquad & \text{min. NO and CO emissions} \\
        \st & \text{Torque}_t \geq \text{PrescribedTorque}_t, \; \forall t=1,\ldots,T \qquad & \text{trip profile} \\
          & \left( \text{NO}_t, \text{CO}_t, \text{Torque}_t \right) = \text{DNN}\left( f_t, r_t, c \right) & \text{DNN constraints} \\
          & f_t \in [f_{\min}, f_{\max}], \; r_t \in [r_{\min}, r_{\max}], \; c \in [c_{\min}, c_{\max}] & \text{bounds on controls}.
    \end{array}
\end{equation}

We are using  1,500 time intervals (corresponding to $T=1500$ seconds in a 25-minute drive). As written, we will have $2T$ + 1 continuous control variables in addition to the $3T$ output variables predicted by the DNN for a total of $5T+1$ continuous variables in this formulation.

Our model has a separate DNN evaluation for each time interval $t$, meaning each evaluation engenders different neural network activations. Because of the presence of the design variable $c$, which is independent of $t$, this problem does not decompose into individual time steps. One may construct a bilevel optimization problem wherein on the upper level we decide $c$ and other engine-level variables of interest and on the lower level we determine the drive-specific variables that change over time by solving $T$ separate optimization problems of much smaller size. This approach will not be addressed in this paper.






\subsubsection{Convex Hull Constraints}

In addition to the constraints that encode the evaluation of a neural network,  constraints must be added to ensure that the solver does not extrapolate significantly from the training data. 
Without these constraints, the optimization routine may find that the optimal solution resides in an area for which the neural network has not learned the behavior of the modeled function, leading to a solution that, while optimal, does not reflect the true function behavior and may be nonsensical. 
In fact, in an earlier implementation without the convex hull constraint, we observed that the optimal design 
was obtained for an engine that produced negative emissions. 
We rectify this error by introducing constraints that constrain the input data to our model to be within the convex hull of the input training data. 

\begin{figure}[htbp]
    \centering
    \includegraphics[width=0.45\linewidth]{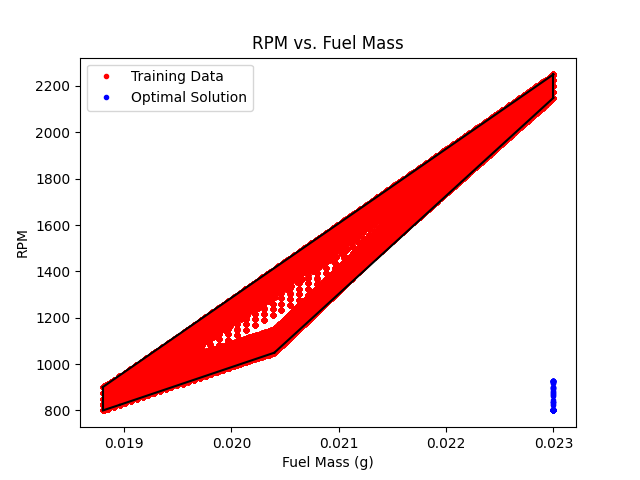}
    \includegraphics[width=0.45\linewidth]{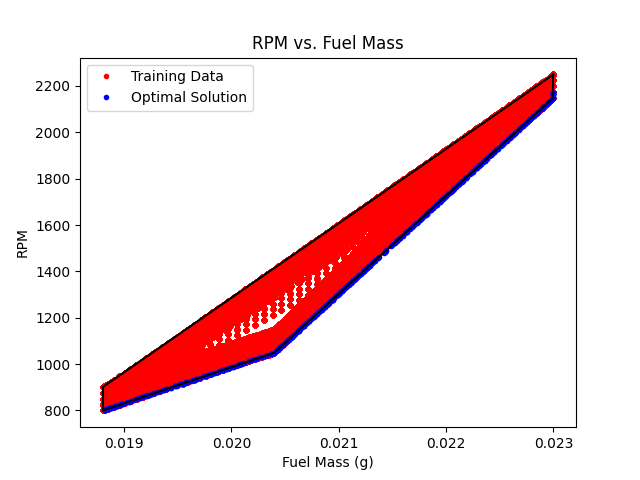}
     \caption{Location of training data and optimal solution for an engine design problem on a neural network with 1 hidden layer of 16 nodes. On the left the input is constrained by box constraints and on the right by convex hull constraints.}
    \label{fig:convex-hull}
\end{figure}

\begin{figure}[htbp]
    \centering
    \includegraphics[width=0.45\linewidth]{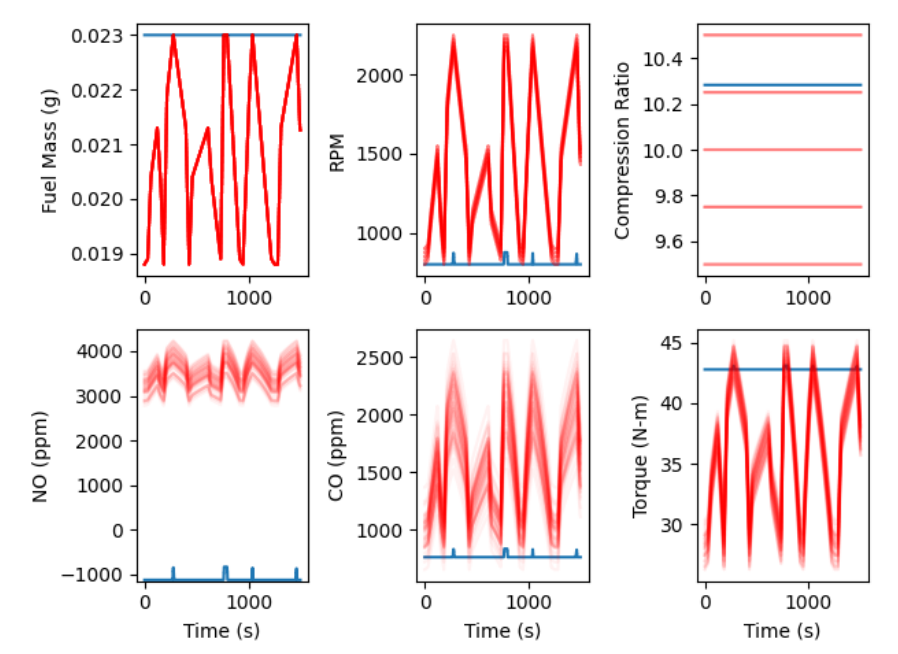}
    \includegraphics[width=0.45\linewidth]{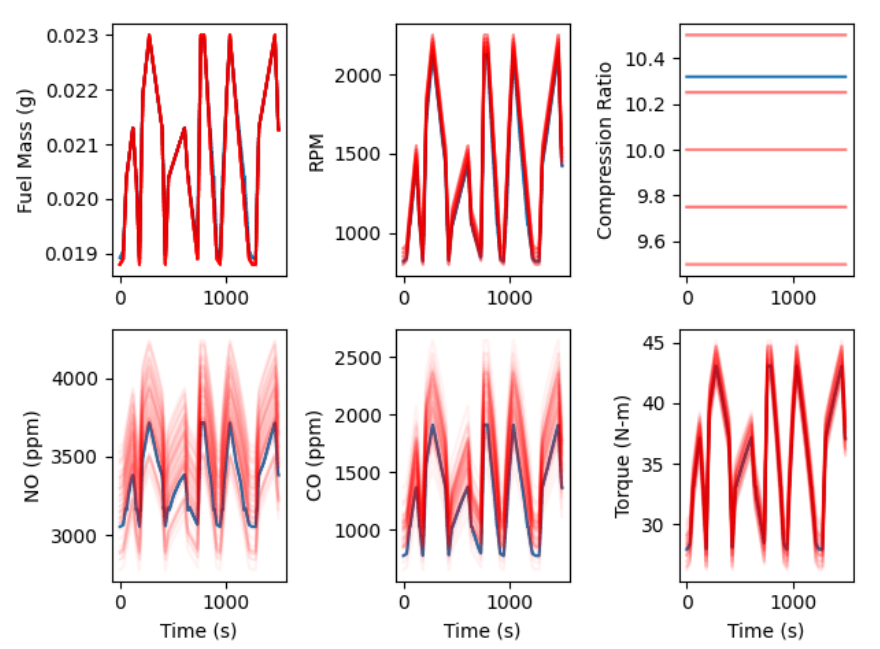}
    \caption{Optimal configuration for each time step when using box constraints (left) and using convex hull constraints (right). The red lines indicate the drive cycles used to train the neural network, and the blue lines indicate the optimal solution.}
    \label{fig:convex-hull-solutions}
\end{figure}

In our experiments we examined how the optimization models operated with simple box constraints that bounded the input by the extremal values of the training data, as well as with the convex hull constraints. Figure~\ref{fig:convex-hull} shows the location of fuel mass and RPM for solutions computed for a sample neural network trained on our data set compared with the actual training data. The corresponding solutions are displayed in Figure~\ref{fig:convex-hull-solutions}.

Without the convex hull constraints, the optimal solution of neural network evaluation appears to wander off to an area for which  no training data exists. 
This situation immediately creates problems in the computed solutions, which are obviously nonsensical because they involve negative emissions. This suggests that our neural network does not generalize well to data that has not yet been seen by the network. 
On the other hand, constraining the input data to be within the convex hull has the effect of producing solutions that look like the training data.

Unfortunately, these constraints have the effect of restricting possible solutions to those that are close to the training data, when the function being modeled may have a minimum that appears far from the training data. 
This situation speaks to the limitations of using solely a surrogate model for optimization. To find a global optimum of the original function would necessarily involve an alternating approach wherein the surrogate model was optimized and then this solution was queried against the modeled function for new data to be added to the surrogate. 
This type of approach is addressed in \citep{queipo2005surrogate} and a similar approach is implemented in \cite{papalexopoulos2021constrained}.

\begin{figure}[htbp]
    \centering
    \includegraphics[width=0.5\textwidth]{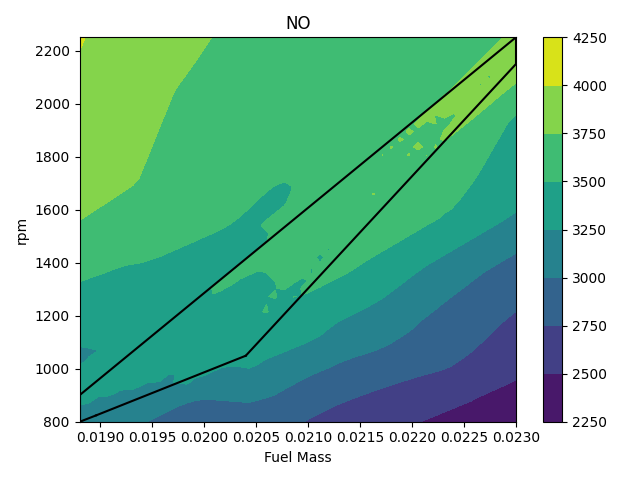}
    \caption{Contour plot for the NO output of a 5-layer DNN.}
    \label{fig:contour}
\end{figure}

If we fix the compression ratio, we can plot the contours of the objective function
as a function of fuel mass and RPM. Figure~\ref{fig:contour} shows the contours of the objective function for a
5-layer DNN. We observe that this objective function is highly nonconvex in the reduced space of the original 
variables but it is convex in the lifted MINLP space according to Proposition~\ref{prop:ConvexDNN}. 
We also 
observe that the convex hull constraint fulfills a second function by acting as a mild convexifier of the 
problem by restricting the variables to a small sliver of the feasible set.

\subsubsection{Warmstarting the Problem}

Seeding a MILP or NLP solver with a feasible solution of high quality can  significantly reduce the solve time. 
For MILP solvers, a feasible solution with sufficiently low value can be used to prune subproblems with  higher objective value in a branch-and-bound tree search. 
For NLP solvers, beginning with a feasible solution avoids the need to search for a feasible solution, and the choice of a good solution can ensure quicker convergence to better locally optimal solutions.

For the optimal engine design problem, we can use our collection of training data to find parameterizations of the engine that exceed the desired torque value but also have a low amount of emissions produced. We do not directly work with the output data from the training set but rather with the output from the neural network applied to the input data since this is what is constrained.

To produce a high-quality solution is then a matter of fixing the parameters that must remain constant for the entire drive (i.e., the compression ratio) and then, from the training data with these fixed parameters, choosing the remaining controls so as to minimize emissions while still exceeding the desired torque for each time step. The remaining variable values correspond to activation of each neuron in the neural network when applied to the input data and hence must be set to exactly those values. This procedure is summarized in Algorithm~\ref{alg:warmstart}.

\begin{algorithm}
    \caption{Algorithm for finding a warmstart solution to \eqref{eq:optEngine}}
    \label{alg:warmstart}
    \KwData{Input $\text{DNN}$, compression ratio $c_0$, training data $\mathcal{X}$, torque profile $\text{TorqueProfile}_t$, $T$}
    Select subset of training data $\tilde{\mathcal{X}}$ with compression ratio $c_0$\;
    Let $\tilde{\mathcal{Y}} = \text{DNN}(\tilde{\mathcal{X}})$\;
    \For{$t=1,\ldots,T$}{
        Let $\text{NO}, \text{CO}, \text{Torque} = \tilde{\mathcal{Y}}$\;
        Let $\mathcal{T}$ be the set of times $t$ where $\text{Torque}_t \ge \text{TorqueProfile}_t$\;
        Let $\tau = \text{argmin}_{t\in\mathcal{T}} (NO_t + \lambda CO_t)$\;
        Let $x_t = \tilde{\mathcal{X}}_\tau$\;
    }
    Set remaining variable values from activation of DNN when applied to $x_t$ for all $t$.\;
\end{algorithm}

\subsection{Adversarial Attack Generation}

Our next problem involves the resilience of neural networks to incorrect classifications.
Deep neural networks have the immensely useful property of being able to uniformly approximate any function in a certain general class, but at the cost of being fairly opaque in terms of how the underlying machinery works.
This opacity may mask unpredictable behavior that makes these neural networks susceptible to attacks that disrupt the correct classifications observed on training and testing data.
Szegedy et  al.~\citep{szegedy2013intriguing} first demonstrated that almost imperceptible perturbations of image data may lead to misclassifications, in essence demonstrating key instabilities in neural networks. This work has led to a number of papers \citep{carlini2016evaluating, cheng2017maximum, fischetti2018deep} developing algorithms that produce estimates of the resilience of neural networks, estimating how close  an incorrectly classified image can be to a correctly classified one.


The central problem of interest as introduced by \citep{szegedy2013intriguing} is described as follows. Given a classifier $\text{DNN}$, an image $x \in [0,1]^m$, and a desired classification label $l$, we have the following problem:
\begin{equation}\label{eq:adversarialAttack}
\begin{array}{lll}
    \dps \mini_{z}\; & \|x - z\| & \\
    \st\; & \text{DNN}(z) =  l & \\
        & z \in [0,1]^m,
\end{array}
\end{equation}
where $\|\cdot\|$ is a suitable norm.
Essentially, we ask for an image of minimal distance to image $x$ that is classified with a different label. 

We add additional constraints (see, e.g., \citep{fischetti2018deep}) to ensure that this is a very confident classification. 
Specifically, we replace the constraint $\text{DNN}(z) = l$ with a constraint asking that the activation for label $l$ is some factor larger than the activation for each of the other labels; in other words, if $y_L$ is the final layer, we have $y_{L,l} \ge \alpha y_{L,i}$ for $i \ne l$. 
Often  the final layer is given by applying a softmax function, namely,
\begin{equation}
y_{L,i} = \sigma(z_1,\ldots,z_n)_i = \exp(z_i)/\sum_{j=1}^n\exp(z_j),
\end{equation}
where the produced values $y_{L,i}$ effectively represent probabilities that the image is a given label $i$.
In this case we can represent this constraint in terms of $z_i$ as
\[
    \sigma(z)_l/\sigma(z)_i \ge \alpha \Leftrightarrow \exp(z_l - z_i) \ge \alpha \Leftrightarrow z_l \ge z_i + \log(\alpha),
\]
which is a linear constraint.
Then if we write $y_L = \text{DNN}(z)$ as the final layer prior to the softmax layer, we will have the following problem:
\begin{equation}\label{eq:adversarialAttack2}
\begin{array}{lll}
    \dps \mini_{z}\; & \|x - z\| & \\
    \st\; & y_L = \text{DNN}(z) & \\
        & y_{L,l} \ge y_{L,i} + \log(\alpha), \qquad i \ne l \\
        & z \in [0,1]^m .
\end{array}
\end{equation}

We observe that the formulation in \eqref{eq:adversarialAttack2} is easily extendable and that constraints on the allowable perturbations can augment the formulation. For example, we may restrict the magnitude of the perturbation to any given pixel by restricting $|x_i - z_i| \le \epsilon$ for all $i$. Alternatively, we may be interested in continuous perturbations and therefore restrict $|(x_i - z_i) - (x_j - z_j)| < \epsilon$ for all pixels $j$ adjacent to pixel $i$. Constraints of this sort have been used to extend adversarial attack optimization problems in \cite{fischetti2018deep}.

\subsection{Surrogate Modeling of Oil Well Networks}

The next example of an optimization problem with an embedded neural network involves the operation of an offshore oil platform and is taken from \citep{grimstad2019relu}. The full optimization problem is reproduced in Eqns.~(\ref{eq:oil-well-problem}), and the associated sets are given in Table \ref{tab:oil-well-sets}.


\begin{table}
    \centering
    \begin{tabular}{c|l} \toprule
        Set & Description \\ \midrule 
        $N$ & Set of  nodes in the network. \\
        $N^w$ & Subset of well nodes in network. \\
        $N^m$ & Subset of manifold nodes in the network. \\
        $N^s$ & Subset of separator nodes in the network. \\
        $E$ & Set of edges in the network. \\
        $E^d$ & Subset of discrete edges that can be turned on or off. \\
        $E^r$ & Subset of riser edges. \\
        $E_i^{\text{in}}$ & Subset of edges entering node $i$. \\
        $E_i^{\text{out}}$ & Subset of edges leaving node $i$. \\
        $C$ & Oil, gas, and water. \\
        \bottomrule
    \end{tabular}
    \caption{Sets used in oil well optimization problem}
    \label{tab:oil-well-sets}
\end{table}

\begin{figure}
    \centering
    \includegraphics[width=0.7\linewidth]{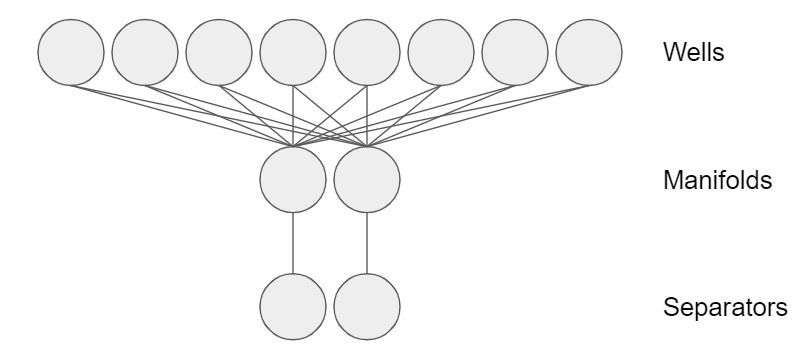}
    \caption{Example oil well network with 8 wells, 2 manifolds, and 2 separators, from \cite{grimstad2019relu}.}
    \label{fig:oil-well}
\end{figure}

In this problem we have a network comprising three sets of nodes: the wells $N^w$ that are the sources of oil, water, and gas; the manifolds $N^m$ that are connected to wells and mix incoming flows of oil, water, and gas; and the separators $N^s$ that are the sinks of all the flow. Each well is connected to each manifold by a pipeline in $E^d$ that may be turned on and off. Each manifold is then connected to a unique separator by a riser in $E^r$.  An example network with 8 wells, 2 manifolds, and 2 separators taken from \cite{grimstad2019relu} is depicted in Figure~\ref{fig:oil-well}.

\begin{subequations} \label{eq:oil-well-problem}
\begin{align}
    \dps \maxi_{y,q,p,\Delta p} & \dps \sum_{e \in E^r} q_{e,\text{oil}} &  \\
    \st & \sum_{e \in E_i^{\text{in}}} q_{e,c} = \sum_{e\in E_i^{\text{out}}} q_{e,c}, \; \forall c \in C, i \in N^m \label{eq:flow-balance}\\
      & p_j = \text{DNN}_e(q_{e,\text{oil}}, q_{e,\text{gas}}, q_{e,\text{wat}}, p_i), & \forall e = (i,j) \in E^r \label{eq:p-nn}\\
      & \Delta p_{e} = p_i - p_j,  & \forall e = (i,j) \in E^r \\
      & - M_e (1 - y_e) \le p_i - p_j - \Delta p_e, & \forall e = (i,j) \in E^d \label{eq:on-off1}\\
      & p_i - p_j - \Delta p_e \le M_e (1 - y_e), & \forall e = (i,j) \in E^d \label{eq:on-off2}\\
      & \sum_{e \in E_i^{out}} y_e \le 1, & \forall i \in N^w \label{eq:routing}\\
      & y_e q^L_{e,c} \le q_{e,c} \le y_e q_{e,c}^U & \forall c \in C, e \in E^d \label{eq:flow-off}\\
      & p_i^L \le p_i \le p_i^U, & \forall i \in N \label{eq:pressure-bound} \\
      & \sum_{e \in E_i^{\text{out}}} q_{e,\text{oil}} = \text{DNN}_i(p_i), & \forall i \in N^w \label{eq:q-nn}\\
      & \sum_{e \in E_i^{\text{out}}} q_{e,\text{gas}} = c_{e,\text{gor}} \sum_{e \in E_i^{\text{out}}} q_{e,\text{oil}}, & \forall i \in N^w \label{eq:ratio1}\\
      & \sum_{e \in E_i^{\text{out}}} q_{e,\text{wat}} = c_{e,\text{wor}} \sum_{e \in E_i^{\text{out}}} q_{e,\text{oil}}, & \forall i \in N^w \label{eq:ratio2}\\
      & p_i = p_i^s = \text{const.}, & \forall i \in N^w  \\
      & y_e \in \{0, 1\}, \forall e \in E^d
\end{align}
\end{subequations}

Our goal in this problem is then to optimize the flow rate of oil to each of the sinks.  Various physical and logical constraints are used to encode the operation of this flow network.  \eqref{eq:flow-balance} ensures that the flow into each node matches the flow out of each node. \eqref{eq:on-off1} and \eqref{eq:on-off2} ensure that if a pipeline is open, then the difference in pressure between the two nodes is actually represented by $\Delta p_e$.  \eqref{eq:routing} ensures that each well  routes its flow only to a single manifold. \eqref{eq:flow-off} bounds the flow rate of each material when pipelines are active and forces the rate to zero when inactive, whereas \eqref{eq:pressure-bound} bounds the pressure at each node. \eqref{eq:ratio1} and \eqref{eq:ratio2} establish the expected ratio of flow rates of each material.

The neural networks $\text{DNN}_e$ and $\text{DNN}_i$ appear in the problem in \eqref{eq:p-nn} and \eqref{eq:q-nn} and represent nonlinear functions that predict the separator pressure and outgoing oil flow rate based on incoming flow rate and pressure, respectively. There is one neural network $\text{DNN}_e$ for each riser edge $e$ and one neural network $\text{DNN}_i$ for each well $i$. Each is trained separately; and for our particular configuration of the network, we have 2 risers and 8 wells for a total of 10 separate neural networks that are encoded into our problem.

The authors in \cite{grimstad2019relu} considered two different configurations of neural networks for this problem: a ``shallow'' network with few layers of many nodes and a ``deep'' network with many layers of fewer nodes. The shallow configuration had two hidden layers of 20 nodes for the well networks $f_i$ and two hidden layers of 50 nodes for the riser networks $g_e$. The deep configurations had four hidden layers of 10 nodes for the well networks and 5 layers of 20 nodes for the riser networks. All in all, this corresponds to 520 ReLU nodes for each configuration.

\section{Embedded Neural Network Formulation}

Given that many deep learning libraries (e.g.,, TensorFlow \cite{abadi2016tensorflow} and PyTorch \cite{paszke2017automatic}) have well-developed built-in automatic differentiation capabilities, we naturally want to see whether we can directly embed the evaluation of the neural network into an NLP:
\begin{equation}\label{eq:embedded-formulation}
    \mini_x \; f(\text{DNN}(x),x) \quad \st \; c(\text{DNN}(x),x) = 0, \; x \in X.
\end{equation}
Unlike the formulations in Section~\ref{sec:form}, this formulation has the advantage of not requiring auxiliary variables for each of the internal nodes of the neural network. That is, this formulation should scale significantly better as the number of nodes in the neural network increases.

The modeling package JuMP \cite{DunningHuchetteLubin2017} in the programming language Julia is a library that establishes a general framework for representing generic optimization problems and interfacing with solvers. 
Of particular use in this problem is its ability to handle nonlinear functions with user-provided gradients. 

The current release of JuMP  supports only univariate user-defined nonlinear functions, but we can represent vector-valued functions by listing each output component separately. A single nonlinear constraint representing $y_i = \text{DNN}(x)_i$ can be written as
\begin{verbatim}
    register(model, DNN_i, n, DNN_i, DNN_i_prime)))
    @NLconstraint(model, y[i] == DNN_i(x...)))),
\end{verbatim}
where \texttt{model} is the JuMP optimization model, \texttt{n} is the dimension of the input, and \texttt{DNN\_i} and \texttt{DNN\_i\_prime} represent function evaluations of the neural network's $i$th output neuron and its gradient that can be provided by TensorFlow.

For networks with  many outputs, explicitly listing these commands  can quickly become unwieldy, but we can enumerate these constraints using macros. 
The following code demonstrates how one can encode $y = M(x)$, where $x \in \mathbb{R}^n$ and $y \in \mathbb{R}^m$, and \texttt{f} and \texttt{f\_p} are evaluation of the neural network and the Jacobian, respectively.
\begin{verbatim}
macro DNNConstraints_grad(model, x, y, n, m, f, f_p)
  ex = Expr(:block)
  for i = 1:m
  dnn = gensym("DNN")
  push!(ex.args, :($dnn_s(z...) = $f(z...)[$i]))
    p = gensym("DNN_prime")
    push!(ex.args,:($prime(g, z...) = begin g .= $f_p(z...)[$i,:] end))
    push!(ex.args,:(register($model, $(QuoteNode(dnn_s)), $n, $dnn, $p)))
    push!(ex.args,:(@NLconstraint($model, $y[$i] == $dnn($x...))))
  end
  ex
end
\end{verbatim}
This macro can then be called to add constraints on a neural network that takes in values from $\mathbb{R}^{20}$ and outputs values in $\mathbb{R}^{10}$ using the following command,
\begin{verbatim}
@DNNConstraints_grad(model, x, y, 20, 10, DNN, DNN_prime),
\end{verbatim}
where \texttt{DNN} and \texttt{DNN\_prime} are now vector- and matrix-valued functions that return the output vector and Jacobian for neural network evaluation, respectively.

With this macro, we now have the ability to directly embed evaluation of a neural network and its derivatives within a mathematical program. 
This allows us to treat neural network evaluation as simply another function that appears in functions and constraints so that we can use any NLP solver for our optimization problems.

\subsection{Convergence Behavior}\label{sec:convergence-behavior}

We have observed in our experiments that state-of-the-art solvers express real difficulties with convergence when applied to this formulation.
We believe that the reason has to do with the choice of activation function, $\text{ReLU(a)} = \max(a, 0)$.
This activation function results in nonconvex constraints and objective in \eqref{eq:OptiML}; moreover, the function is nonsmooth whenever $a = 0$.

The nonsmoothness of the ReLU function has generally been considered to be of little concern since, in practice, neural networks involving ReLU neurons have found great success especially in problems of classification \cite{he2016deep,simonyan2014very,krizhevsky2012imagenet} so that they have become an industry standard in deep learning. Part of their success is attributed to their ease in implementation as well as their tendency to produce sparse activation patterns and avoid the vanishing gradient problem in training \cite{glorot2011deep}. Furthermore,  some theoretical results have affirmed some convergence guarantees under mild conditions \cite{li2017convergence, du2018gradient} and \cite{goodfellow2014qualitatively} have demonstrated qualitatively that the training problem often experiences few of the issues common in nonconvex and nonsmooth optimization.

\begin{figure}
    \centering
    \includegraphics[width=0.45\linewidth]{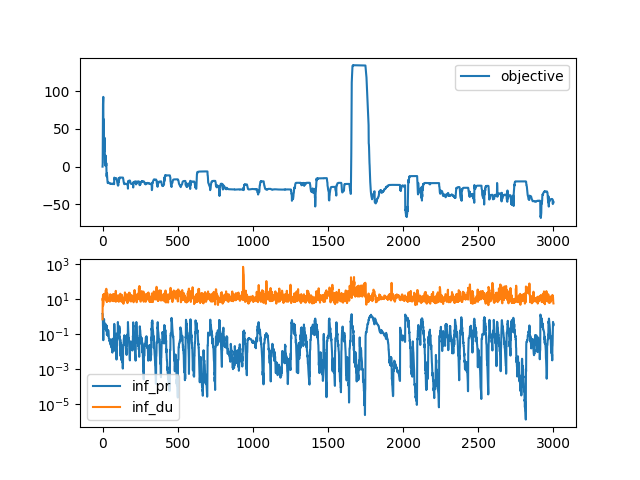}
    \includegraphics[width=0.45\linewidth]{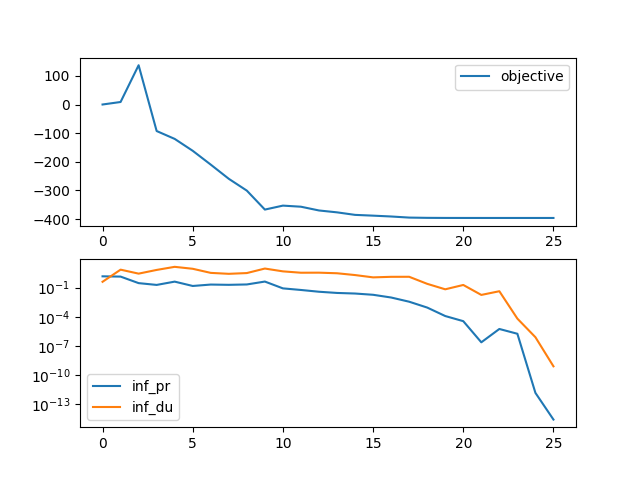}
    \caption{Objective value, primal infeasibility, and dual infeasibility plotted against iteration number for a sample ReLU network (left) and a sample swish network (right) over an Ipopt solve.}
    \label{fig:dual-infeas}
\end{figure}

These results pertain primarily to the optimization problems solved in the training process. We have observed relatively little work in the literature regarding optimization problems involving trained networks. 
Our preliminary experiments suggest that, in this context, the expectation of good behavior may be unfounded. 
Initial solves of Problem (\ref{eq:optEngine}) using the state-of-the-art NLP solver Ipopt \cite{wachter.biegler:06} experienced serious issues with convergence. 
Figure \ref{fig:dual-infeas} illustrates the objective value and primal and dual infeasibility for one example solve of the problem. Qualitatively, we observe that instead of terminating, the solver bounces about among objective values of about the same magnitude.
The reason for the lack of convergence can be explained by the fact that one measure of convergence, the dual infeasibility, remains high throughout the duration of the solve. 
This type of behavior makes this formulation difficult for general use because there are not clear conditions in general for checking whether the current solution at hand is ``good.'' 
We suspect that these convergence issues can arise in part due to the nondifferentiability of the ReLU neurons. 
Classical optimization algorithms such as quasi-Newton algorithms \cite{powell1969method} often have basked in the assumption of differentiability, and convergence results typically depend on them.

An alternative type of activation function that has seen some success in classifier networks \cite{ramachandran2017swish} is the swish function, which is defined as follows:
\begin{equation}
    \text{swish}(a) = \frac{a}{1 + e^{-\beta a}}. \label{eq:swish-def}
\end{equation}
$\beta$ is a hyperparameter that may be learned but is typically set to 1. Note that for $\beta = 0$ this is simply a linear function and that for $\beta \to\infty$ this approaches the ReLU function. With the activation function in \eqref{eq:swish-def}, all functions in \eqref{eq:OptiML} are twice continuously differentiable under Assumption \ref{ass:general}.

To ensure that failure to converge when using the ReLU network is not due to choice of NLP solver, we perform the same solve of Problem (\ref{eq:optEngine}) using instead swish neurons. We  plot the objective and infeasibilities in Figure \ref{fig:dual-infeas}. 
Instead of failing to terminate after 3,000 iterations, Ipopt converges successfully after only 25 iterations; and both primal infeasibility and dual infeasibility converge rapidly to zero.
As a differentiable network, the swish formulation does not express the convergence issues that the nondifferentiable ReLU networks appear to exhibit.

These initial  results confirm the observation in \cite[Fig. 1]{bolte2020conservative}
that showed that ReLU networks can be {\em nondifferentiable almost everywhere} for certain
types of networks. This observation makes the use of ReLU networks questionable as
embedded constraints within optimization solvers that rely on differentiable problem functions. Hence, in Section \ref{sec:form} we
develop alternative formulations that have better numerical properties.

\subsection{Stationarity in the Embedded Formulation}

In this section we will give a condition for the stationarity of the embedded formulation as well as prove a specific form for the Jacobian of a ReLU neural network that holds under mild genericity conditions.
To handle the points of nondifferentiability in the embedded formulation, we consider generalized Jacobians first introduced in \cite{clarke.et.al:96}, which can be defined for all $x$ in the domain of a function $g$ as follows:
\begin{equation}\label{eq:clarke-grad-def}
    \partial g(x) = \text{conv}\{ \lim_k \nabla g(x_k) : x_k \rightarrow x, x_k \notin S \},
\end{equation}
where $S$ is the set of points of nondifferentiability of $g$ and $\text{conv}(A)$ is the convex hull of the set $A$. Since ReLU networks are piecewise linear functions by construction, 
the generalized Jacobian will be given by a convex hull of the Jacobians of linear functions on finitely many regions in the neighborhood of a given point.
Stationarity of $g$ at a point $x$ is then equivalent to having $0 \in \partial g(x)$.

\begin{figure}
    \centering
    \includegraphics[width=\linewidth]{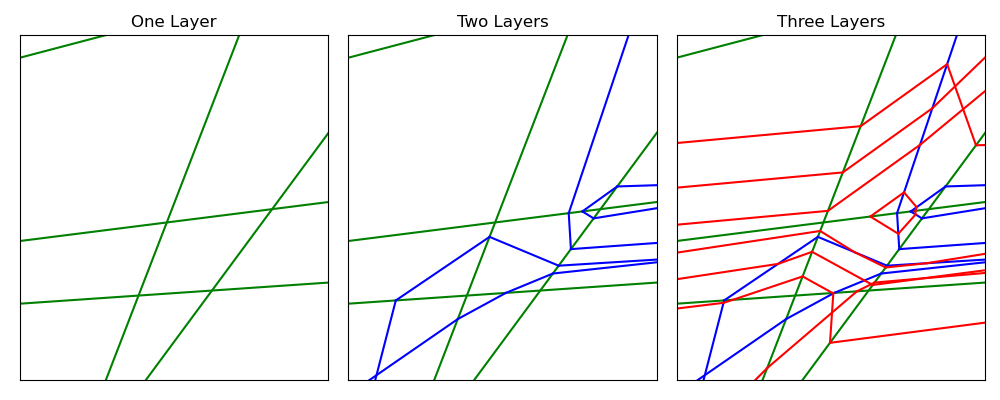}
    \caption{Regions of a randomly instantiated neural network with two inputs followed by 1, 2, and 3 layers of 5 ReLU neurons each.}
    \label{fig:regions}
\end{figure}

We elaborate on what we mean by region in the context of feedforward ReLU networks. 
As these neural networks are created by composing sequences of ReLU activation functions with affine maps, we can specify the region of the function by deciding which ReLU neurons will be active; each selection corresponds to a different affine function. 
Given a set of active neurons $N$ (we will refer to these as \textit{activation patterns}), we will define the associated region $R(N)$ to be the set of points $x$ such that the value at a neuron in $N$ prior to activation is positive and this value for any other neuron is negative. The value prior to activation being zero corresponds to being on the boundary of another region.
$R(N)$ will be empty if no point $x$ engenders the activation pattern $N$.
The regions will tile the input space, and adjacent regions are separated by the 0-level sets of a particular neuron.

Depictions of the regions of neural networks with 1, 2, and 3 layers can be seen in Figure \ref{fig:regions}.
We observe that the regions of a single-layer network arise from arrangements of hyperplanes but thst as more layers are added, the regions become more complex as they are subdivided.
The total amount of nonempty regions of both shallow and deep neural networks has been studied in detail in \cite{montufar2014number} and \cite{pascanu2013number}. However, their definition of region as maximally connected linear regions of the piecewise function slightly differs from ours because we distinguish adjacent regions that have the different activation patterns but ultimately get mapped to the same linear map.

Before we can discuss stationarity of the embedded formulation, we will need to introduce some notation. We will identify any neuron by its layer index $\ell$ and the index $i$ within the layer.
Given a set of neurons in the neural network $N = \{(\ell_1, i_1), \ldots, (\ell_n, i_n)\}$, we define $F_{N}$ to be the affine form given when the set of neurons in $N$ are active and all other neurons are inactive.
Put precisely, given the affine transform in layer $\ell$, $A^{\ell}x = W^{\ell^T}x + b^{\ell}$, if neuron $(\ell, i)$ is inactive, we replace column $w_i$ and component $b_i$ with zeros. We define our new affine form as
\[
    \hat{A}^{\ell}_N(x) = \hat{W}^{\ell^T}x + \hat{b}^{\ell}, \; \hat{W}^{\ell} = W\text{diag}(\kappa^\ell), \; \hat{b} = \text{diag}(\kappa^\ell)b.\,
\]
where $\kappa^{\ell}$ is a vector with $\kappa^{\ell}_i$ equal to 1 if neuron $(\ell, i) \in B$, and zero otherwise. Then $F_N(x) = \hat{A}^L_N\hat{A}^{L-1}_N\cdots\hat{A}^1_Nx$. 

We can partition the neurons of layer $\ell$ into three sets given $x$: strictly inactive neurons $I^{\ell,-}(x)$, nonstrictly inactive neurons $I^{\ell,0}(x),$ and active neurons $I^{\ell,+}(x)$. 
We define these to be the indices $i$ of the neurons in layer $\ell$ for which $A^{\ell}_i\hat{A}^{\ell-1}\cdots\hat{A}^1x$ is negative, zero, and positive, respectively.
The sets $I^+(x), I^0(x), I^-(x)$ are then defined to be all layer-neuron pairs $(\ell, i)$, where $i \in I^{\ell,+}(x), I^{\ell,0}(x),$ or $I^{\ell,-}(x)$, respectively.
With these definitions we observe that for any $x$, $\text{DNN}(x) = F_{I^+(x)}(x)$.

Note that $F_N$ is now smooth as a composition of affine functions, and we can compute the Jacobian $J_{F_N}$ simply by the chain rule:
\begin{equation}
    J_{F_N} = \hat{W}^1\hat{W}^2\cdots\hat{W}^L.
\end{equation}

In light of \eqref{eq:clarke-grad-def}, the generalized gradient at $x$ can then be written as the convex hull of these Jacobians associated with the activation pattern for each region in the neighborhood of $x$. 
Any such activation pattern must necessarily include $I^+(x)$ and exclude $I^-(x)$ but may only include a subset of $I^0(x)$.
If there is a region corresponding to each set of the form $I^+(x) \cup N$, where $N \subset I^0(x)$, we can write the generalized gradient in a nice form given by the following proposition.
\begin{proposition} \label{prop:convex-hull}
Suppose at a point $x$ that the region $R(I^+(x) \cup M)$ is nonempty for each choice of $M \subset I^0(x)$. Let $S$ be the set of all matrices of the form
\[
    W^1\text{diag}(\kappa^1)\cdots W^L\text{diag}(\kappa^L),
\]
where
\[
    \begin{cases}
        \kappa^{\ell}_i = 1 & (\ell, i) \in I^+(x), \\
        \kappa^{\ell}_i \in [0,1] & (\ell, i) \in I^0(x), \\
        \kappa^{\ell}_i = 0 & (\ell, i) \in I^-(x).
    \end{cases}
\]
Then the generalized gradient of $\text{DNN}(x)$ at $x$, $\partial \text{DNN} (x) = S$.
\end{proposition}
\begin{proof}
First, we remark that any Jacobian $J_{F_N}$, where $N = I^+(x) \cup M$, is clearly in $S$ since we can take $\kappa^{\ell}_i = 1$ for all $(\ell, i) \in M$ and take it equal to 0 for $(\ell, i) \in I^0(x) \setminus M$, and the two expressions agree. Furthermore any convex combination of the $J_{F_B}$ should be contained in $S$ since $S$ is a convex set. Hence $\partial \text{DNN}(x) \subset S$.

Now given a matrix $V = W^1\text{diag}(\kappa^1)\cdots W^L\text{diag}(\kappa^L)$ in $S$, we show that it is in $\partial \text{DNN}(x)$. We enumerate $I^0(x) = \{(\ell_1, i_1), \ldots, (\ell_n, i_n)\}$. 
Consider first any alteration $V_1$ of $V$ by fixing $\kappa^{\ell_1}_{i_1}$ and replacing $\kappa^{\ell_j}_{i_j}$ by $\hat{\kappa}^{\ell_j} \in \{0,1\}$ for $j > 1$. 
Note that $V_1$ is a convex combination of the two matrices given by replacing $\kappa^{\ell_1}_{i_1}$ with 1 and by replacing it with 0. 
Both these matrices are in $\partial \text{DNN}(x)$ since this set includes all Jacobians associated with any activation pattern on $I^0(x)$. Hence we must have that $V_1 \in \partial \text{DNN}(x)$. If we then consider fixing the first two components and replacing the remaining such components to produce $V_2$, in the same fashion we can find two matrices of the form $V_1$ of which $V_2$ is a convex combination and hence $V_2 \in \partial \text{DNN}(x)$. Repeating until we have fixed all components, we will have that $V \in \partial \text{DNN}(x)$.
\end{proof}


Unless the region is nonempty for each of these activation patterns, the generalized Jacobian will not necessarily take this form.
As a counterexample, consider a zero-bias two-layer network with one input variable and $W^1 = (1, 1)$, $W^2 = (1, -1)^T$. 
This corresponds to having $\text{DNN}(x) = \text{ReLU}(x) - \text{ReLU}(x) = 0$, and so the generalized gradient is $\{0\}$ everywhere, whereas the formulation given above would introduce many more possible gradients.

\begin{figure}
    \centering
    \includegraphics[width=0.5\linewidth]{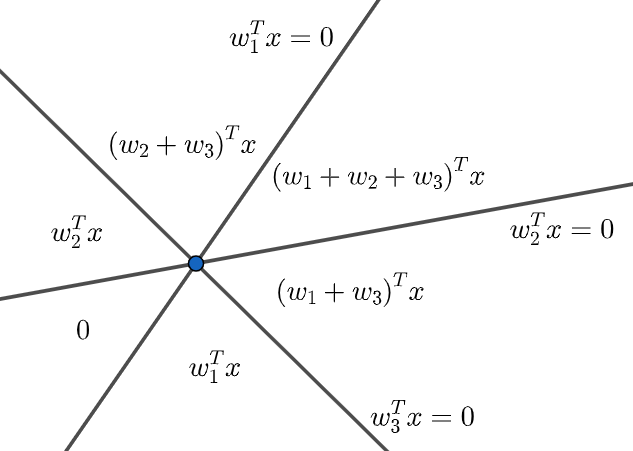}
    \caption{Regions of a zero-bias, single-layer network with two input variables and three output variables.}
    \label{fig:three-neuron}
\end{figure}

For a slightly more interesting example, consider a one-layer network with two input variables, no bias, and three output variables, namely, $\text{DNN}(x) = \text{ReLU}(W^Tx) = \text{ReLU}(w_1^Tx) + \text{ReLU}(w_2^Tx) + \text{ReLU}(w_3^Tx)$. The input domain $\mathbb{R}^2$ is then divided into regions by the lines $w_i^Tx = 0$ for $i=1,2,3$ each intersecting at the origin, as can be seen in Figure \ref{fig:three-neuron}. 
Observe that in the neighborhood of the origin, there are only six neighboring regions and  two regions are missing, one where the output would be $(w_1 + w_2)^Tx$ and another where the output is $w_3^T$, corresponding to activation patterns where only the first two neurons are active and only the third neuron is active, respectively. In this scenario the generalized gradient is similarly not as expansive as we desire.

Requiring each region $R(I^+(x)\cup N)$ to be nonempty seems to be an onerous requirement, especially given there are $2^{|I^0(x)|}$ such regions. However, it turns out that under some mild genericity conditions, we can expect this to be true.
These conditions arise from the theory of hyperplane arrangements. 
Suppose we have a collection of hyperplanes $H_1,\ldots,H_n \subset \mathbb{R}^d$, where $H_i$ is given by $\{x \in \mathbb{R}^d: a_i^Tx = b_i\}$ for some vectors $a_i \in \mathbb{R}^d$ and constants $b_i \in \mathbb{R}$. 
We say the collection is in \textit{general position} if the intersection of $m$ hyperplanes has dimension $d - m$ if $m \le d$ and is empty if $m > d$. If the intersection of all hyperplanes is nonempty, it is called a \textit{central arrangement}.
In particular, if $d = 2$, the hyperplanes are lines in $\mathbb{R}^2$, and they are in general position if no three lines intersect at a point and no two lines are parallel. 
A well-known result from Zaslavsky \cite{zaslavsky1975facing} is that if we have $m$ hyperplanes in general position, the number of regions they divide $\mathbb{R}^d$ into is given by $\sum_{i=0}^d \binom{m}{i}$. Note that if $m \le d$, this quantity is $2^m$, and we can identify each region by choosing whether $a_i^Tx < b$ or $a_i^Tx > b$ for all $i$.

To ensure that the necessary regions are not empty, we will require a specific set of hyperplanes to be in general position. Before we describe the hyperplanes, we will prove a useful lemma.

\begin{lemma} \label{lemma:general-position}
Suppose a central arrangement of hyperplanes $H_1,\ldots, H_n$, where $H_i = \{x \in \mathbb{R}^n: a_i^Tx = b_i\}$. Then the set $H_1,\ldots,H_{n-1},H_n^A$, where $H_n^A = \{x: (\sum_{i\in A} w_ia_i + a_n)^Tx = \sum_{i \in A}w_ib_i + b_n\}$ for any choice of constants $w_i \in \mathbb{R}$, is also in general position for any $A \subset \{1,\ldots, n-1\}$.
\end{lemma}
\begin{proof}
Because the hyperplanes intersect, we can change coordinates so that each hyperplane intersects the origin; thus,  without loss of generality, we can write $b_i = 0$ for all $i$. In this setting, hyperplanes being in general position is equivalent to the normal vectors being linearly independent. Obviously if $a_1,\ldots, a_n$ are linearly independent, $a_1,\ldots,a_{n-1},a_n + \sum_{i\in A}w_i a_i$ are linearly independent.
\end{proof}

Using this lemma, in the following proposition we can now show exactly which hyperplanes  need to be in general position.

\begin{proposition} \label{prop:general-position}
Suppose that the set of hyperplanes $\{ H_{\ell,i} = \{x: A^{\ell}_i\hat{A}^{\ell-1}_{I^+(x)}\cdots \hat{A}^{1}_{I^+(x)} x = 0\}: \ell=1,\ldots, L, i \in I^{0,\ell}(x)\}$ is in general position. Then for any choice of $N \subset I^0(x)$, the set of hyperplanes $\{ H_{\ell,i} = \{x: A^{\ell}_i\hat{A}^{\ell-1}_{M}\cdots \hat{A}^{1}_{M} x = 0\}: \ell=1,\ldots, L, i \in I^{0,\ell}(x)\}$ is also in general position with $M = I^+(x) \cup N$.
\end{proposition}
\begin{proof}
Observe for any hyperplane $A^{\ell}_i\hat{A}^{\ell-1}_{I^+(x)}\cdots \hat{A}^{1}_{I^+(x)} x = 0$ and any $j$ in $1, \ldots, \ell-1$, we consider the associated affine form $A^{\ell}_i\hat{A}^{\ell-1}_{I^+(x)}\cdots\hat{A}^{j+1}_{I^+(x)}(\hat{A}^j_{I^+(x)} + \hat{A}^{j}_{N})\hat{A}^{j-1}_{I^+(x)}\cdots \hat{A}^{1}_{I^+(x)} x = 0$. Let $b^Tx + c$ denote the affine form associated with $A^{\ell}_i\hat{A}^{\ell-1}_{I^+(x)}\cdots\hat{A}^{j+1}_{I^+(x)}$, and for each neuron $i$ in layer $j$ let $d_i^Tx + e_i$ denote the affine form associated with $\hat{A}^j_{i}\hat{A}^{j-1}_{I^+(x)}\cdots \hat{A}^{1}_{I^+(x)} x$. Then we can write
\[
A^{\ell}_i\hat{A}^{\ell-1}_{I^+(x)}\cdots\hat{A}^{j+1}_{I^+(x)}(\hat{A}^j_{I^+(x)} + \hat{A}^{j}_{N})\hat{A}^{j-1}_{I^+(x)}\cdots \hat{A}^{1}_{I^+(x)} x = \sum_{i \in I^{+,j}(x) \cup N_j}b_i (d_i^Tx + e_i) + c_i,
\]
and by assumption the hyperplane defined by $H = \{x: \sum_{i \in I^{+,j}(x)}b_i(d_i^Tx + e_i) + c_i = 0\}$ is in general position with the hyperplanes defined by $H_i = \{x: d_i^Tx + e_i = 0\}$ for $i \in N_j$. The new hyperplane is created by adding to the normal vector of $H$ a linear combination of the normal vectors of each $H_i$, and so it   remains in general position with the $H_i$. Hence, we have replaced one $I^+(x)$ with $M$. Repeating in this fashion, we can replace each $I^+(x)$ in each affine form with $M$ and still remain in general position.
\end{proof}

Now that for any choice of $N \subset I^0(x)$ the associated hyperplanes will be in general position, we can be assured that the associated region is nonempty.

\begin{proposition}\label{prop:regions}
Suppose the hypothesis of Proposition \ref{prop:general-position} holds at a given $x$. Then for each subset $N \subset I^0(x)$, the region $R(M)$ with $M = I^+(x)\cup N$ is nonempty.
\end{proposition}
\begin{proof}
From Proposition \ref{prop:general-position}, the hyperplanes $\{ H_{\ell,i} = \{x: A^{\ell}_i\hat{A}^{\ell-1}_{M}\cdots \hat{A}^{1}_{M} x = 0\} : \ell=1,\ldots, L, i \in I^{0,\ell}(x)\}$ are in general position; and since these form a central arrangement (centered at $x$), there must be a nonempty region with $A^{\ell}_i\hat{A}^{\ell-1}_{M}\cdots \hat{A}^{1}_{M} x > 0$ if $(\ell, i) \in M$ and $A^{\ell}_i\hat{A}^{\ell-1}_{M}\cdots \hat{A}^{1}_{M} x < 0$ otherwise. 
\end{proof}

The condition that the hyperplanes enumerated in Proposition \ref{prop:general-position} be in general position is fairly mild.
Certain types of neural networks will violate this assumption: for  example, zero-bias neural networks will violate this assumption since each hyperplane will go through the origin.
Additionally, we cannot have $|I^0(x)|$ be larger than the dimension of the input space since this also will have too many intersecting hyperplanes.
That being said, we generally expect a collection of randomly generated hyperplanes to be in general position with probability 1 under most natural probability distributions.

With these tools in hand, we can now state in full the stationarity conditions for the embedded problem:

\begin{theorem}
Suppose at a given $x^*$, the conditions of Proposition \ref{prop:general-position} hold. Then $x^*$ is a stationary point of \eqref{eq:embedded-formulation} if  there exist nonnegative multipliers $\mu^*$ and matrices $\hat{W}^1,\ldots,\hat{W}^L$ with $\hat{W}^1\cdots\hat{W}^L \in \partial \text{DNN}(x)$ such that
\begin{enumerate}
    \item $c(\text{DNN}(x^*),x^*) \le 0$
    \item $\mu^{*^T}c(\text{DNN}(x^*),x^*) = 0$
    \item $\nabla_{x}f^* + \nabla_{x}c^*\mu^* + \hat{W}^1\cdots\hat{W}^L(\nabla_{y}f^* + \nabla_{y}c^*\mu^*) = 0$,
\end{enumerate}
where $f^*, c^*$ indicate evaluation of $f, c$ at $x^*$.
\end{theorem}
\begin{proof}
The first two conditions are standard primal feasibility and complementary slackness conditions. To arrive at the third, we note that the gradient of the Lagrangian $L(x, \mu) = f(\text{DNN}(x), x) + \mu c(\text{DNN}(x), x)$ on each region (with activation pattern $B$) nearby is given by $\nabla_x f + \nabla_x c\mu + J_{F_B}(\nabla_y f + \nabla_y c\mu)$ evaluated at $x$. 
As $x\to x^*$ in a given region, and since $f$ and $c$ are smooth, this approaches $\nabla_x f^* + \nabla_x c^*\mu + J_{F_B}(\nabla_y f^* + \nabla_y c^*\mu)$. 
As the generalized gradient at $x^*$ is given by the convex hull of all such gradients and by Proposition \ref{prop:regions} there is one for region $I^+(x)\cup A$ with $A \subset I^0(x)$, it follows in a similar fashion as in the proof of Proposition \ref{prop:convex-hull} that any element of the generalized gradient can be written $\nabla_{x}f^* + \nabla_{x}c^*\mu^* + \hat{W}^1\cdots\hat{W}^L(\nabla_{y}f^* + \nabla_{y}c^*\mu^*)$. 
Thus condition 3 holds.
\end{proof}
\section{Formulating DNNs as Optimization Models}\label{sec:form}

We provide two alternative formulations of ReLU DNNs in terms of optimization models
that avoid the pitfalls of the embedded formulation in Section \ref{sec:convergence-behavior}. The first
formulation uses binary variables to model the $\max$-functions, resulting in a
mixed-integer program, building on \cite{grimstad2019relu,fischetti2018deep,anderson2019strong}.
Our formulation differs from \cite{cheon2020outer}, which  considered 
only inverse problems and did not
exploit the convex structure of the DNNs that arises when we lift the DNN constraint, resulting in a 
nonconvex mode. The second formulation uses complementarity constraints
that can be solved as systems of nonlinear inequalities. We derive theoretical
properties of both formulations.

\subsection{Formulating DNNs with Mixed-Integer Sets}\label{sec:formMIP}

In this section we show how general optimization problems involving DNNs, such as \eqref{eq:OptiML},
can be equivalently formulated as convex MIPs, extending \cite{fischetti2018deep}. 



We assume that the DNN is a deep neural network with ReLU activation functions, and we rewrite the
nonconvex problem \eqref{eq:OptiML} as a constrained problem:
\begin{equation}\label{eq:ConstrLossReLU}
    \mini_{x,y^L} \; f(x,y^L) \quad \st \; y^L = \text{DNN}(x), \; c(x,y^L) \leq 0, \; x \in X,
\end{equation}
where $L$ is the number of layers of the DNN.
Fischetti  and Jo~\citep{fischetti2018deep}  have shown
that the nonconvex constraint, $y^L = \text{DNN}(x)$, can be formulated as a mixed-integer linear
set in the case of DNNs with ReLU activation functions. We let $w^{\ell}_{i}$
denote the weights of neuron $i=1,\ldots,N_{\ell}$ at level $\ell=1,\ldots,L$ and $b^{\ell}_{i}$ its corresponding
bias. Then, the levels $\ell=1,\ldots,L$ are computed as
\begin{equation}\label{eq:ReLU}
    y^{\ell}_{i} = \text{ReLU}\left(w^{\ell^T}_{i} y^{\ell-1} + b^{\ell}_{i} \right)
        = \max\left(0, w^{\ell^T}_{i} y^{\ell-1} + b^{\ell}_{i} \right).
\end{equation} 
We lift this constraint by introducing slack variables, $s^{\ell}_{i}$, and binary variables, $z^{\ell}_{i} \in
\{0,1\}$, and observe that \eqref{eq:ReLU} is equivalent to the
mixed-integer linear constraints
\begin{equation}\label{eq:ReLU-MIP}
    y^{\ell}_{i} - s^{\ell}_{i} = w^{\ell^T}_{i} y^{\ell-1} + b^{\ell}_{i}, \;
    0 \leq y^{\ell}_{i} \leq M^{\ell,y}_{i}(1-z^{\ell}_{i}), \;
    0 \leq s^{\ell}_{i} \leq M^{\ell,s}_{i} z^{\ell}_{i}, \; z^{\ell}_{i} \in \{0,1\},
\end{equation} 
where $M^{\ell,y}_{i}, M^{\ell,s}_{i}>0$ are sufficiently large constants (if $z^{\ell}_{i}=1$, we are on the $0$-branch of
ReLU, and if $z^{\ell}_{i}=0$, we are on the positive branch). 
By substituting \eqref{eq:ReLU-MIP} into \eqref{eq:ConstrLossReLU} we obtain a convex MINLP
that is equivalent to the problem \eqref{eq:OptiML}:
\begin{equation}\label{eq:OptiDNN}
\begin{array}{lll}
    \dps \mini_{x,y,s,z}\; & f(x,y_L) & \\
    \st\; & c(x,y_L) \leq 0 & \\
          & y^{\ell}_{i} - s^{\ell}_{k} = w^{\ell^T}_{i}y^{\ell-1} + b^{\ell}_{i}, \; & \ell=1,\ldots,L \\
        & 0 \le y^{\ell}_{i} \le M^{\ell,y}_{i}(1 - z^{\ell}_{i}), \; 0 \le s^{\ell}_{i} \le M^{\ell,s}_{i}z^{\ell}_{i}, \; &  i=1,\ldots,N_{\ell}, \ell=1,\ldots,L \\
        & y_0 = x, \; z^{\ell}_{i} \in \{0,1\}, x \in X. &
\end{array}
\end{equation}
We observe that we have as many binary variables in this problem as we have ReLU nodes. 
Next we show that the resulting lifted formulation results in a tractable convex MINLP.

\begin{proposition}\label{prop:ConvexDNN}
Let Assumption~\ref{ass:general} hold. Then it follows that \eqref{eq:OptiDNN} is a convex MINLP in the sense
that the continuous relaxation of \eqref{eq:OptiDNN} is a convex NLP.
\end{proposition}

\begin{proof} The result follows from the convexity of $f$, $c$, and $X$ and the fact that the remaining
constraints are affine (with the exception of the integrality restriction on $z$).
\end{proof}
    
The choice of the big-$M$ constants $M^{\ell}_{i}$ in (\ref{eq:ReLU-MIP}) can have a great impact on the solution time of our MIP. If $M^{\ell}_{i}$ is too small, the problem excludes solutions that should be feasible; but if it is too large, the space that must be searched by the solver may become so large as to be too computationally intractable for most computers.

These upper bounds are common in mixed-integer programs, and we follow the approach in \citep{fischetti2018deep} where we consider each neuron in our neural network individually, removing all constraints on other neurons in the same or subsequent layers. Then we set as our objective function to maximize $y^{\ell}_{i}$ in one iteration and $s^{\ell}_{i}$ in another iteration. These computed optimal values will serve as upper bounds $M^{\ell,y}_{i}$ and $M^{\ell,s}_{i}$, respectively, and can be used in our model. 
This process of obtaining bounds is related to the optimality-based bound-tightening technique used in 
global optimization; see, for example, \cite{gleixner2017three}.
Since these constants depend solely on the weights for the neural network as well as inputs to the neural network, we  need  to compute these bounds only once and may reuse them in any optimization problem involving this neural network. Alternative methods for handling the big-$M$ constraints are considered in \citep{grimstad2019relu}, and an alternative formulation with exponentially many constraints alleviating the above concerns is presented in \citep{anderson2019strong}. 

With this formulation established, we may then easily pass the problem as is to any standard MINLP solver such as Gurobi \citep{gurobi5}, CPLEX \citep{cplex12}, Bonmin \cite{bonami2007bonmin}, MINOTAUR \cite{MINOTAUR}, or Baron \cite{sahinidis:96} to compute the solution. Because MINLPs are NP-complete, however, these problems do not scale well, and only problems involving modestly sized networks (on the order of 100 hidden nodes) may be tractably solved. Neural networks used in commercial settings generally involve at least thousands of hidden nodes, resulting in thousands of binary variables, a number that  typically is well beyond the reach of commercial solvers.

\subsection{Formulating DNNs with Complementarity Constraints}

In this section we discuss an alternative formulation of \eqref{eq:OptiML} as a nonconvex nonlinear program  using complementarity constraints. This approach has the advantage of scaling significantly better for problems with larger neural networks, but with the caveat that solutions produced can  be guaranteed only to be locally optimal. 

Our approaches are based on the following observations.
We can rewrite the ReLU activation function in \eqref{eq:ReLU}
equivalently as a complementarity constraint (using vector notation):
\begin{equation}\label{eq:ReLU-CC}
        y^{\ell} = \max( W^{\ell^T} y^{\ell-1} + b^{\ell} , 0 )
        \; \Leftrightarrow \;
        0 \leq y^{\ell} \; \perp \; y^{\ell} \geq W^{\ell^T} y^{\ell-1} + b^{\ell},
\end{equation} 
where $\perp$ means that for each component, $i$, both inequalities $y^{\ell}_{i}\ge 0$ and   $y^{\ell}_{i} \ge [W^{\ell^T} y^{\ell-1} + b^{\ell}]_i$ are satisfied and at least one is satisfied at equality. By replacing the ReLU function with these 
complementarity constraints, we obtain a mathematical program with complementarity
constraints, which we can solve using standard NLP solvers; see, for example, \citep{Leyf03a,LeyLopNoc:06,FLRS:06,RaghBieg:05}. The mathematical program modeling language {\tt AMPL} \citep{fourer.gay.kernighan:93} allows the modeling of complementarity constraints: for example, the above constraint can be written as follows:
\begin{verbatim}
ReLUCompl{l in Level, i in Neuron[l]}: 0 <= y[l,i] complements 
y[l,i] >= sum{j in Neuron[l-1]} W[l,i,j] y[l-1,j] + b[l,i] .
\end{verbatim}
    
The most successful NLP solvers handle MPCCs by reformulating the complementarity constraints in
\eqref{eq:ReLU-CC} by first introducing the same slack variables as for the MINLP and lifting 
the formulation
    \[
    s^{\ell} = y^{\ell} - w^{\ell^T}_{i} y^{\ell-1} + b^{\ell},
    \]
    and then rewriting \eqref{eq:ReLU-CC} equivalently as
    \[
    s^{\ell} = y^{\ell} - w^{\ell^T}_{i} y^{\ell-1} + b^{\ell}\; \text{and} \; 0 \leq y^{\ell} \; \perp \; s^{\ell} \geq 0.
    \]
    A nonlinear optimization formulation is then given as
    \[
    s^{\ell} = y^{\ell} - w^{\ell^T}_{i} y^{\ell-1} + b^{\ell}, \quad y^{\ell} \geq 0,\quad  s^{\ell} \geq 0, \; \text{and} \;
    y^{\ell^T} s^{\ell} \leq 0,
    \]
    where the lower bound, $y^{\ell^T} s^{\ell} \geq 0$, is implied by the nonnegative bounds on $y^{\ell}, s^{\ell} \geq 0$,
    and again omitted for numerical reasons; see \citep{FLRS:06}.

Because the last constraint involves the nonconvex term $y^{\ell^T}s^{\ell}$, the problem as a whole is nonconvex, and standard NLP solvers will produce only locally optimal solutions instead of globally optimal ones.
Altogether, we then have the following NLP:
\begin{equation}\label{eq:OptiDNN-NLP}
\begin{array}{lll}
    \dps \mini_{x,y,s,z}\; & f(x,y_L) & \\
    \st\; & c(x,y_L) \leq 0 & \\
          & y^{\ell}_{i} - s^{\ell}_i = w^{\ell^T}_{i}y^{\ell-1} + b^{\ell}{i}, \; & \ell=1,\ldots,L \\
        & 0 \le y^{\ell}_{i}, \; 0 \le s^{\ell}_{i}, \; &  i=1,\ldots,N_{\ell}, \ell=1,\ldots,L \\
        & \dps\sum_{\ell=1}^{L}\sum_{i=1}^{N_{\ell}} y^{\ell}_{i}s^{\ell}_{i} \le 0,  \\
        & y_0 = x, \; x \in X. &
\end{array}
\end{equation}
Note that we could have used the constraint $y^{\ell^T}s^{\ell} \le 0$ for each $\ell$ separately. We prefer the formulation in \eqref{eq:OptiDNN-NLP} because it has better convergence behavior in practice; see \cite{FletLeyf:04}.
Again we have a continuous variable for every node in our neural network as well as a slack variable for each ReLU node. The absence of integer variables, however, produces a much more scalable problem.

We now demonstrate that the stationarity conditions for the MPCC formulation coincide with those of the embedded formulation. In MPCC form, the optimization problem is given by the following:
\begin{equation}\label{eq:OptiMultiMPCC}
    \begin{array}{ll}
     \dps \mini_x \; & f\left(x^0,x^L)\right) \\
     \st       \; & c\left(x^0,x^L)\right) \leq 0 \\
                & 0 \leq x^{\ell} \; \perp \; x^{\ell} \geq W^{\ell^T} x^{\ell-1} + b^{\ell}, \; \ell=1,\ldots,L.
    \end{array}
\end{equation}

First, we state the definition of strong stationarity for the MPCCs  \eqref{eq:OptiMultiMPCC}; see, for  example, \cite{SchSch:00} for its general form.  We then show that the two conditions are equivalent in our case.

\begin{definition}[Scheel and Scholtes,  \cite{SchSch:00}]\label{def:strongStat}
We say that $(x^{0,*},x^{1,*},\ldots,x^{L,*})$ is a strongly stationary point of \eqref{eq:OptiMultiMPCC} if there exist multipliers $\mu^* \geq 0$ and $\nu_1^*, \nu_2^*$ such that the following conditions are satisfied:
\begin{subequations}\label{eq:StrongStatMulti}
 \begin{align}
     & c^* \leq 0 \; \text{and} \; 0 \leq x^{\ell,*} \; \perp \; x_{\ell,*} \geq W^{\ell^T} x^{\ell-1,*} + b^{\ell}, \; \ell=1,\ldots,L \\
     & \mu^{*,T} c^* = 0 \\
     & 0 = \nabla_{x^0} f^* + \nabla_{x^0} c^* \mu^* + W^1 \nu_2^{1,*} \label{eq:KKTx1} \\
     & -\nu_1^{\ell,*} - \nu_2^{\ell,*} + W^{\ell+1}\nu_2^{\ell+1,*} = 0, \; \ell=1,\ldots,L-1 \label{eq:KKTx2} \\
     & \nabla_{x^L} f^* + \nabla_{x^L} c^*\mu^* - \nu_1^{L,*} - \nu_2^{L,*} = 0 \label{eq:KKTx3} \\
     & x^{\ell,*}_j > 0 \Rightarrow \nu^{\ell,*}_{1j} = 0 \; \text{and} \; x^{\ell,*}_j > w^{\ell^T}_j x^{\ell-1,*} + b^{\ell}_j \Rightarrow \nu^{\ell,*}_{2j} = 0, \; \ell=1,\ldots,L\\
     & 0 = x^{\ell,*}_j = w^{\ell,T}_j x^{\ell-1,*} + b_j \Rightarrow \nu^{\ell,*}_{1j} \geq 0, \nu_{2j}^* \geq 0, \; \ell = 1,\ldots,L , \label{eq:biactiveMulti}
 \end{align}   
\end{subequations}
where all functions and gradients are evaluated at $(x^*,y^*)$, that is,   $f^* := f(x^*,y^*)$.
\end{definition}

Using these stationarity conditions, we can arrive at the following conditions on the gradients of our optimization problem.

\begin{proposition}\label{prop:strong-stationarity}
Given a strongly stationary point satisfying \eqref{eq:StrongStatMulti}, there exist $\hat{W}^1,\ldots,\hat{W}^L$, column scalings of $W^1,\ldots,W^L$ such that $\nabla_{x^0}f^* + \nabla_{x^0}c^*\mu^* + \hat{W}^1\cdots\hat{W}^L(\nabla_{x^L}f^* + \nabla_{x^L}c^*\mu^*) = 0$.
\end{proposition}
\begin{proof}
To construct such matrices, we work backward starting with layer $L$ and examine each component $j = 1,\ldots,N_L$. Backsubstituting $\nu^{L,*}_{2}$ using \eqref{eq:KKTx2} and \eqref{eq:KKTx3}, we have
\[
    -\nu^{L-1,*}_1 - \nu^{L-1,*}_2 + W^L(\nabla_{x^L}f^* + \nabla_{x^L}c^*\mu^* - \nu^{L,*}) = 0.
\]
For a given component $j$ we can determine the appropriate scaling factor $\kappa^L_j$ in three cases:
\begin{enumerate}
    \item $w^{L^T}_j x^{L-1,*} + b^L_j < 0$: Then it follows that $x^{L,*}_j=0$, $\nu^{L,*}_{2j} = 0$ and $\nu_{1j}^*$ is not sign-constrained. It follows from \eqref{eq:KKTx3} that $\nu^{L,*}_{1j}=[\nabla_{x^L} f^* - \nabla_{x^L} c^* \mu^*]_j$ and hence that $\kappa^L_j=0$.
    \item $w^{L^T}_j x^{L-1,*} + b^{L}_j > 0$: Then it follows that $x^{L,*}_j=w^{L,T}_j x^* + b_j > 0$ and hence that $\nu^{L,*}_{1j}=0$. Thus, we can set $\kappa^L_j=1$.
    \item $w^{L^T}_j x^{L,*} + b^L_j = 0 = x^{L,*}_j$: Then it follows that $\nu^{L,*}_{1j}, \nu^{L,*}_{2j} \geq 0$ from \eqref{eq:biactiveMulti}. From \eqref{eq:KKTx3}, we obtain that $[\nabla_{x^L} f^* + \nabla_{x^L} c^* \mu^*]_j \geq \nu^{L,*}_{1j}, \nu^{L,*}_{2j} \geq 0$; and
    given $\nu_{1j}^* \in \left[0,[\nabla_y f^* + \nabla_y c^* \mu^*]_j\right]$, we choose $\kappa^L_j = \frac{[\nabla_{x^L} f^* + \nabla_{x^L} c^* \mu^*]_j - \nu^{L,*}_j}{[\nabla_{x^L} f^* + \nabla_{x^L} c^* \mu^*]_j}$, which is in $[0,1]$.
\end{enumerate}
Hence we can write $W^L(\nabla_{x^L}f^* + \nabla_{x^L}c^*\mu^* - \nu^{L,*}) = \hat{W}^L(\nabla_{x^L}f^* + \nabla_{x^L}c^*\mu^*)$, where $\hat{W}^L = W\text{diag}(\kappa^L)$. Repeating this process for each $\ell$, we arrive at $\nabla_{x^0}f^* + \nabla_{x^0}c^*\mu^* + \hat{W}^1\cdots\hat{W}^L(\nabla_{x^L}f^* + \nabla_{x^L}c^*\mu^*) = 0$.

\end{proof}

We can now prove the following theorem.
\begin{theorem}\label{theo:EquivStatMulti}
Given $x^*$, suppose the assumptions of Proposition \ref{prop:general-position} hold.
Then $x^*$ is a stationary point of \eqref{eq:embedded-formulation} if and only if $(x^{0,*},x^{1,*},\ldots,x^{L,*})$ is a strongly stationary point of \eqref{eq:OptiMultiMPCC} with $x^{\ell,*}$ being the activations of layer $\ell$ in the neural network.
\end{theorem}

\begin{proof}
Since the given expression lies in the generalized gradient for \eqref{eq:embedded-formulation} and equals zero, strong stationarity of the MPCC formulation clearly implies stationarity in the embedded formulation. 
On the other hand, if $0$ is in the generalized gradient for the embedded formulation, as we know the form of the generalized gradient, we can find $\hat{W}^1,\ldots, \hat{W}^L$ such that
$\nabla_{x^0}f^* + \nabla_{x^0}c^*\mu^* + \hat{W}^1\cdots\hat{W}^L(\nabla_{x^L}f^* + \nabla_{x^L}c^*\mu^*) = 0$. 
The choice of scaling factors $\kappa$ for the matrices will then determine the values of $\nu$ as in Proposition \ref{prop:strong-stationarity} so that the conditions for strong stationarity hold.
\end{proof}

\section{Numerical Experiments}\label{sec:numerics}

Here we present our numerical experiments with neural-network surrogates for the three
sample applications introduced in Section~\ref{sec:model}.

We assess the performance of solvers for the MIP, MPCC, and embedded formulations of \eqref{eq:optEngine}. 
We use the commercial solver CPLEX \citep{cplex12} to solve the MIP formulation of the problem and the solver Ipopt \citep{wachter.biegler:06} to solve the MPCC and embedded formulations of the problem. Both solvers are run with the default options. 
All experiments are performed on a single thread on an Intel Xeon Gold 6130 CPU and 188 GB of memory.

\subsection{Numerical Experiments with Engine Design Optimization}

We have applied our models and algorithms to a collection of neural networks with varying numbers of layers to study how well each formulation scales. 
The architectures we tested all have the simple structure of an input layer with the three input variables, followed by $n$ hidden layers of 16 nodes, and then the output layer with the three output variables. 
We consider networks with 1, 3, and 5 hidden layers and train each of these networks  for 20 epochs on the simulation data produced by the simulator of \citep{aithal2019maltese}. We use the \texttt{adam} solver in TensorFlow \cite{abadi2016tensorflow} to train the neural network.

The number of auxiliary variables for the MIP and MPCC formulations scales with the number of ReLU neurons as well as with the number of time steps. 
The last layer has no ReLU neurons since we want the final layer to be able to take all real values, but each of the hidden layers uses ReLU neurons.
We  then have $T \times (\# \text{Hidden Layers}) \times (\# \text{Nodes Per Layer})$ additional sets of auxiliary variables.

Solving the full integer program with 1,500 time steps to optimality is generally intractable (even for a single-layer network, this amounts to $1500\times16\times1 = 24000$ binary variables, beyond the scope of any general-purpose state-of-the-art solver), so we consider instead a coarser discretization using larger time steps.

Instead of using data at each second as is presented in the original data, we
consider time intervals of 750, 500, 250, 150, 30, 15, 10, 5, and 1 seconds. 
Since we have 1,500 seconds of data, this corresponds to 3, 6, 10, 50, 100,
150, 300, and 1,500 time steps for which we are evaluating our neural network.
At each of the larger time steps, our prescribed torque profile will be the average of the prescribed torques in that interval. 
This approach gives us separate problems for each choice of coarseness of the discretization and for each choice of neural network architecture.
As the number of time steps increases, we have observed that the computed solutions converge to the solution of the full problem.

For each instance, we run until convergence to a solution, 3 hours have passed, or Ipopt has performed 3,000 iterations. 
If the solver has not converged at the end of the experiment, we record the best solution discovered that is feasible and  that we determined to be having a constraint violation under $10^{-6}$.

\begin{figure}
    \centering
    \includegraphics[width=\linewidth]{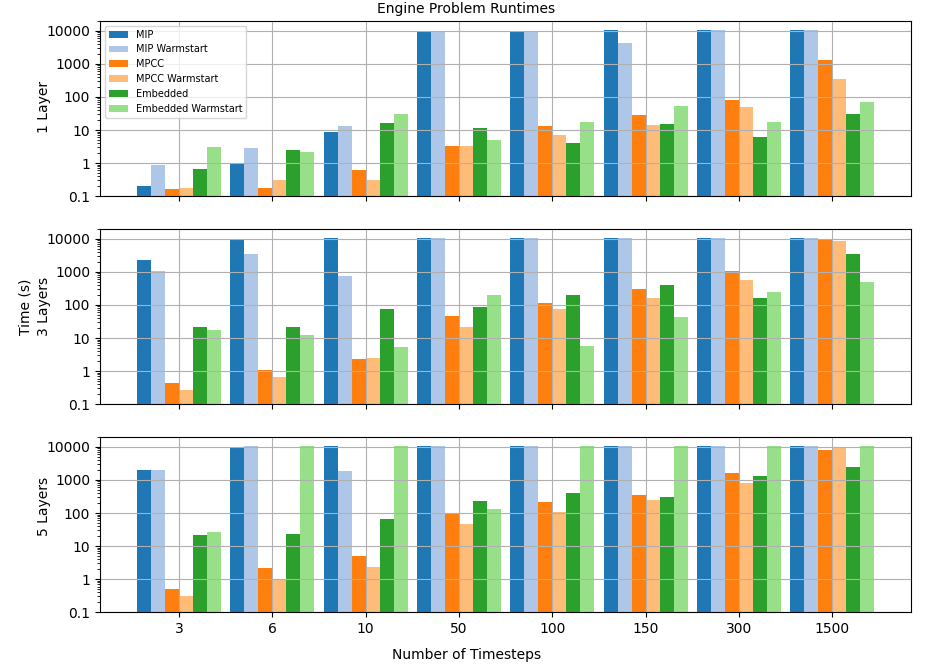}
    \caption{Time until a solver found its best solution within 3 hours (10,800 seconds) on each architecture for varying values of $T$. If a runtime of 10,800 seconds is recorded, this indicates the solver failed to find a feasible solution in the full 3 hours (if a warmstart is used, this means it failed to find a feasible solution better than the warmstart).}
    \label{fig:engine_times}
\end{figure}

\begin{figure}
    \centering
    \includegraphics[width=\linewidth]{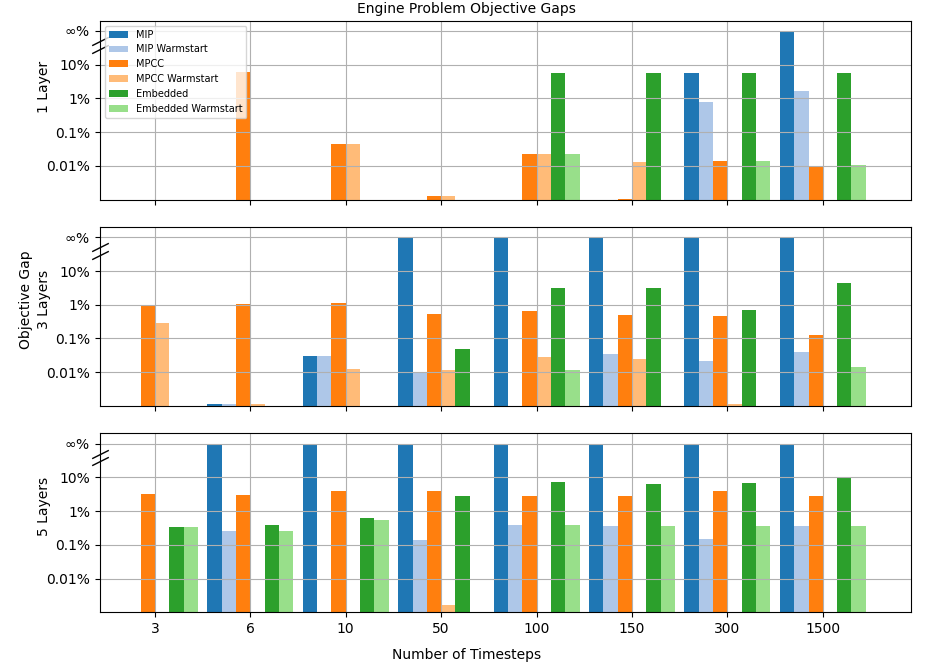}
    \caption{Percentage gap in objective between final objective value and best-known objective value for each solver on each architecture for varying values of $T$. A value of 0 indicates the solver found the best-known objective value, and a value of $\infty$ indicates that no solution was found.}
    \label{fig:engine_gaps}
\end{figure}

Figure~\ref{fig:engine_times} shows the amount of time needed until a solver found its best solution with a time limit of 3 hours on each architecture with each time step size. The percentage gap between the computed objective and the best-known objective for each solver is shown in Figure~\ref{fig:engine_gaps}.

Except for the  smallest problems, CPLEX fails to converge to the optimal solution of the MIP within the allotted time limit of 3 hours.
Even worse, it fails to even find a feasible solution for most problems with architectures with 3 or 5 hidden layers. 
This result is to be expected because, except for the smallest cases, these problems can involve thousands of binary variables. 
With the warmstart, CPLEX is actually able to find solutions with significantly better objective values for each of the configurations. 
It still fails to prove optimality for any of these solutions, however, and times out for the same set of problems as without the warmstart. 

For the MPCC, on the other hand, Ipopt always finds locally optimal solutions within the time limit, often in  seconds for problems with a small number of time steps. 
These solutions may be suboptimal, however, worse by up to 10\% in some instances. With the warmstart solution, Ipopt solves the MPCC in terms of speed generally by a factor of 2. 
Furthermore, the objective solution tends to have a significantly smaller objective, although there are a few exceptions. 
Solving the MIP without warmstart finds a better solution than does  the MPCC formulation in only three problem instances, and in all three
cases the MPCC solution is within 0.1\% of the MIP's globally optimal solution. 

In almost all cases the embedded network formulation is able to find its best solution more quickly than does either the MIP or the MPCC formulation. 
Without the warmstart, however, the embedded ReLU network formulation performs worse in all 
but a few cases, with an objective value that is up to 10\% larger than the MPCC formulation objective values. 
With the warmstart, the objective values are comparable except for the largest
network, where in most cases Ipopt failed to find any solutions better than the warmstart solution in the allotted time, and so the computed objectives are worse by a small margin.
 Ipopt finds the best solution to the embedded formulation  faster than does the MIP formulation, but it is often slower than the MPCC formulation.
Each iteration is faster since the formulations involve significantly fewer variables; 
but because of issues with convergence, the solver often takes significantly more steps. 
This situation is noticeable for the 5-layer network when warmstarting, since Ipopt often fails to converge to any solution better than the warmstart.

Overall, these experiments confirm our prior suspicions that  using the MILP formulation without warmstart quickly becomes computationally intractable as the size of the neural network increases past modest architectures. 
Switching to the MPCC formulation of the problem offers significant speedup at the cost of losing global optimality, although when the MPCC formulation obtains a worse objective value, 
it tends to be by only a marginal amount, and when it does better, the  improvement can be by a significant amount. 
The embedded formulation, on the other hand, can provide speed and scalability but often encounters difficulty with convergence and finds slightly worse objective values. 
For all formulations, significant gains can be realized in terms of both  solving time and solution quality by providing the solver with a high-quality warmstart solution.

\subsection{Numerical Experiments with Adversarial Attack Generation}

Next we consider each of the formulations of the adversarial attack generation problem in \eqref{eq:adversarialAttack2} by considering neural networks trained on the MNIST handwritten digit recognition data set \citep{lecun1998mnist}. 
We consider 10 different architectures corresponding to having an input layer with $28\times28$ nodes followed by either 1 or 2 fully connected layers with 20, 40, 60, 80, or 100 hidden ReLU nodes each and then the output layer with 10 nodes and a softmax activation function. 

If we use the mixed-integer formulation from \eqref{eq:OptiDNN} for the $\text{DNN}$ constraints, then depending on our choice of
norm, it becomes a mixed-integer linear program (for $L^1$ or $L^{\infty}$ norms) or a mixed-integer quadratic program (for $L^2$ norm) that we can solve using CPLEX. 
We will use the $L^2$ norm in experiments. 

For each neural network, we train the model by minimizing the categorical cross-entropy loss function for 10 epochs using the 60,000 digits of training data. 
Each neural network attains a test accuracy of 96--98\% when tested on the 10,000 digits of testing data.

We set $\alpha = 1.2$, meaning the probability for the given classification $l$ is 1.2 times higher than for any other. 
On each architecture, we solve the problem for 100 digits from the training data, with the goal of finding the closest image to the given digit that will be classified as a zero. 
As a warmstart, we initialize the solution in each solve to be a digit that is classified as the desired digit, so that the solver starts at a feasible point.
We run each iteration until convergence to the globally optimal solution (for the MIP formulation) or to a locally optimal solution (for the MPCC and embedded network formulations) or until one hour has passed or 3,000 iterations have occurred in Ipopt, at which point we terminate with the best feasible solution seen so far. 


In Table~\ref{tab:adversarial-results} we tabulate the results of the experiments for all three formulations. 
For each neural network architecture we record the average solve times and objectives for CPLEX and Ipopt over the 100 solves. 
We also record how many times CPLEX finds the optimal solution to the MIP as well as how many times in each formulation a feasible solution that is better than the initial solution is found.
For the embedded formulation we also include the average objective excluding the infeasible problems since we observed a significant difference between these cases.


\begin{table}
\centering
\begin{tabular}{l|cccc|ccc|cccc}
\toprule
 Architecture & \multicolumn{4}{c|}{MIP} & \multicolumn{3}{c|}{MPCC} & \multicolumn{4}{c}{Embedded ReLU}  \\
 &  \pbox{2cm}{Avg. \\ Time} &  \pbox{2cm}{Avg. \\ Obj.} &  \pbox{2cm}{Num. \\ Opt.} &  \pbox{2cm}{Num. \\  Feas.} &  \pbox{2cm}{Avg. \\Time} &  \pbox{2cm}{Avg. \\Obj.}  & \pbox{2cm}{Num. \\  Feas.} & \pbox{2cm}{Avg. \\Time} &  \pbox{2cm}{Avg. \\Obj.} & \pbox{2cm}{Avg. \\ Feas. \\ Obj.}  & \pbox{2cm}{Num. \\  Feas.}\\
\midrule
784,20,10       & 1.8                & \textbf{1.79}            & 100                               & 100                                & \textbf{0.3}                & 1.80             & 100 &   35.5 &       14.67 &     1.80 &            88\\
784,40,10       & 76.9               & \textbf{1.75}            & 100                               & 100                                & \textbf{1.5}                & 1.77       & 100     &   28.4 &       8.64 &    1.77 &             94\\
784,60,10       & 2302             & \textbf{1.85}            & 52                                & 100                                & \textbf{4.8}                & 1.87          & 100  &   31.8 & 1.91     &    1.91 &              100\\
784,80,10       & 2571             & \textbf{2.04}            & 38                                & 100                                & \textbf{14.2}               & 2.07          & 100 &   48.1 &       17.35 &   2.08 &               86\\
784,100,10      & 3078             & \textbf{2.16}            & 16                                & 100                                & \textbf{22.9}               & 2.18           & 100  &   30.1 &       10.66 &  2.09 &                92\\
784,20,20,10   & 106.1              & \textbf{1.20}            & 100                               & 100                                & \textbf{0.8}                & 1.22          & 100 &   26.6 &       4.04 &  4.04 &             100\\
784,40,40,10   & 3100             & 1.90            & 16                                & 100                                & \textbf{4.4}                & \textbf{1.77}           & 100 &   38.8 &       19.97 &      1.68 &            84\\
784,60,60,10   & 3513             & 6.85            & 5                                 & 100                                & \textbf{14.3}               & \textbf{1.88}           & 100 &   29.8 &       12.61 &     12.61 &            100\\
784,80,80,10   & 3587             & 27.59           & 1                                 & 76                                 & 36.1               & \textbf{1.46}     & 100  &   \textbf{27.8} &       3.83 &       3.83 &           100\\
784,100,100,10 & 3600             & 54.83           & 0                             & 42                                 & 57.0               & \textbf{1.78}            & 100         &   \textbf{35.1} &       1.83 &      1.83    &           100\\
\bottomrule
\end{tabular}
\caption{Comprehensive results for solves of the MIP formulation using CPLEX, the MPCC formulation using Ipopt, and the embedded ReLU formulation also using Ipopt. Times are recorded in seconds.}
\label{tab:adversarial-results}
\end{table}

The average time to convergence for the CPLEX solves of the MIP formulation is significantly higher than that of the Ipopt solves of the MPCC formulation. 
Most CPLEX solves for the larger networks time out at 1 hour before finding an optimal solution; for the smaller networks, it can take 10 to 100 times longer to converge to a solution. 
Ipopt, on the other hand, takes less than a minute on average to converge to a locally optimal solution of the MPCC formulation. 
The embedded formulations are fairly uniform in how long they take to find their best solution, which makes sense given that the problem does not change in size with the number of neurons in the network.

The longer solve times of the MIP formulation have the advantage of eventually uncovering better solutions the majority of the time even if they are not provably optimal. 
For all the single hidden layer architectures and the two smallest double hidden layer architectures, the MIP formulation almost always produces a better solution than the MPCC formulation does. 
These solutions generally offer only a marginal improvement on the order of 1--2\% difference in the objective. For perturbations of this magnitude, the differences are essentially imperceptible. 

\begin{figure}
    \centering
    \includegraphics[width=0.5\linewidth]{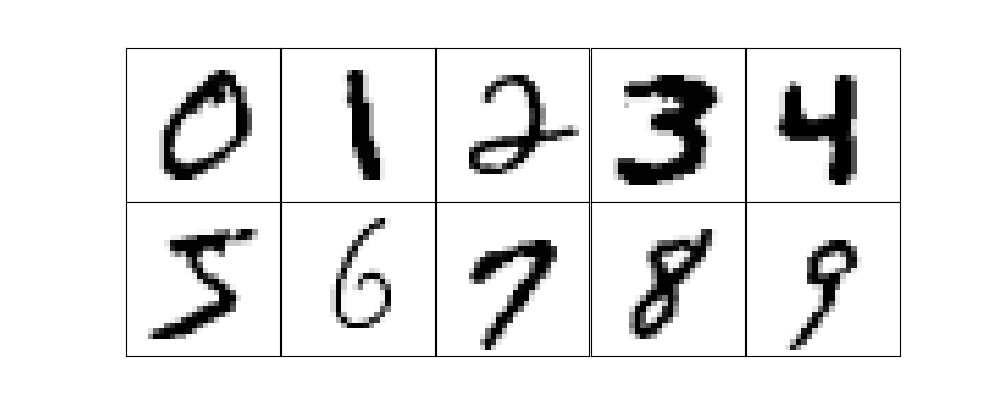}
    \caption{Suite of 10 digits from the MNIST training data for which adversarial attacks are generated.}
    \label{fig:digit_suite}
\end{figure}

\begin{figure}
    \centering
    \includegraphics[width=0.45\linewidth]{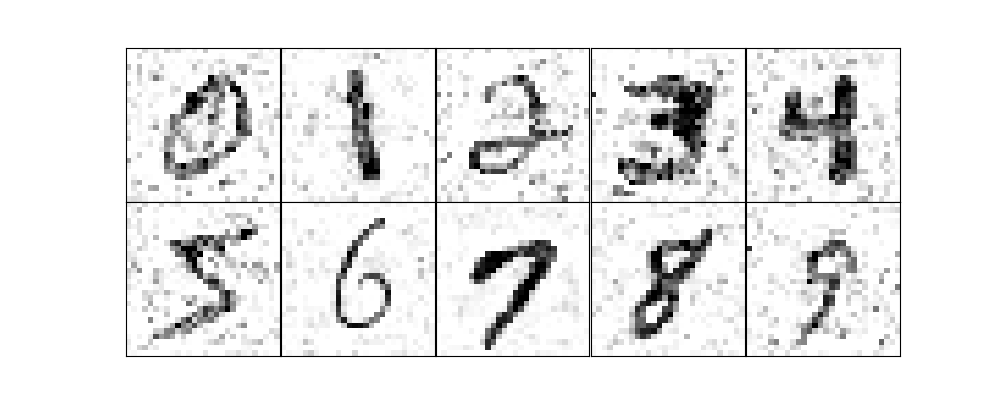}
    \includegraphics[width=0.45\linewidth]{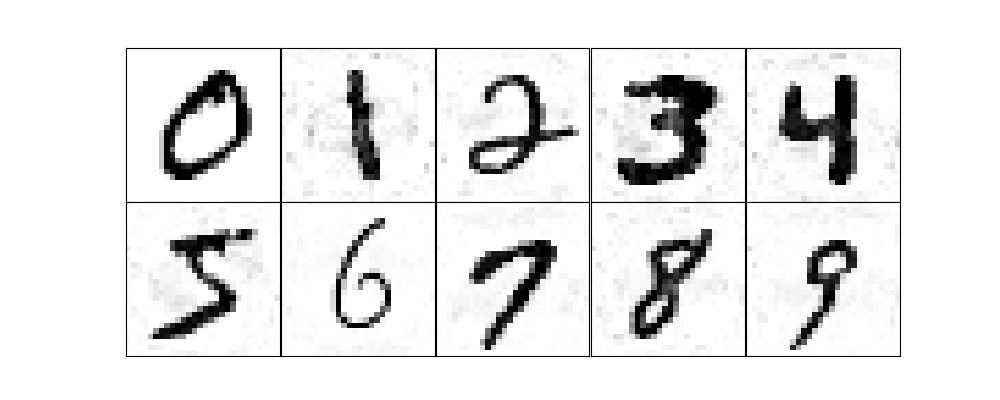}
    \caption{Perturbed images produced after 1 hour of computation when using CPLEX to solve the MIP formulation (left) and Ipopt to solve the NLP formulation (right) with given digits in Figure~\ref{fig:digit_suite} for the neural network with 2 hidden layers of 100 nodes.}
    \label{fig:perturbed_digit_suite}
\end{figure}

For the largest networks the MIP formulation fares poorly and cannot find a better feasible solution than the initial solution in all cases for the two largest networks. 
The feasible solutions it does find are of noticeably lower quality: of all instances from the three largest networks, CPLEX found the best solution in only 13 cases out of 300. 
The MIP solutions for the larger networks are visually worse, as can be seen in a comparison between the perturbed images generated from CPLEX and Ipopt in Figure~\ref{fig:perturbed_digit_suite} with respect to the original images in Figure~\ref{fig:digit_suite}. 

When compared with the MIP and MPCC formulations, the embedded formulation, while often faster and more scalable, can be more unstable in the sense that while using it, Ipopt can fail to find feasible solutions even though for this problem they certainly exist. 
The ReLU networks are more prone to failure and exhibit the poor convergence behavior as in Figure~\ref{fig:dual-infeas} leading to failure to find solutions in about 5\% of instances. In a few cases the objective value of the found solutions is also noticeably worse.

\subsection{Numerical Experiments with Oil Well Networks}

Next we consider solving the oil well problem using each of our formulations.
We again used CPLEX for the MIP formulation and solved the problem with a time limit of 1 hour, recording the best time for each neural network configuration.
This problem involves binary variables even in the complementary constraint formulation, so we needed to use an MINLP solver and elected to use the solver Bonmin \cite{bonami2007bonmin}. We used the default branch-and-bound scheme with Ipopt as the NLP solver.
Because no MINLP solver is currently supported on JuMP for handling both integer variables and user-defined nonlinear functions, we could not solve the problem as stated using the embedded neural network formulations. 

The results for the solves on the full problem are presented in Table \ref{tab:full-oil-well-problem}. We observe that in this small set of problems, the MPCC formulation outperforms the MIP formulation in both network configurations  in terms of time taken as well as objective value produced. CPLEX fails to even produce a solution for the deep neural network configuration. These are hard problems; indeed,  in 3 of the 4 solves performed, the solutions found are not proven optimal.

\begin{table}
    \centering
    \begin{tabular}{c|cc|cc}
    \toprule
     &\multicolumn{2}{c|}{Shallow} & \multicolumn{2}{c}{Deep} \\
    Solver    & Time (s) & Objective & Time (s) & Objective \\ \midrule
    CPLEX &  $3388^*$ & 1.271 & $3600^*$ & NaN \\
    bonmin &  $1865^*$ & \textbf{1.283} & 1275.3 & \textbf{1.304} \\
        \bottomrule
    \end{tabular}
    \caption{Time to find the best solution for each instance of the full oil well problem with a time limit of 1 hour. The * indicates that the solver could not prove the solution optimal.}
    \label{tab:full-oil-well-problem}
\end{table}

To obtain a problem that did not have binary variables, we considered fixing each of the binary variables randomly to either 0 or 1 while still ensuring feasibility. 
If we represent the neural networks using complementary constraints or using the embedded formulation,  the problems become standard NLPs that we can directly solve using Ipopt as before. 
In our setup we considered 10 different configurations of the binary variables, and we solved them using CPLEX for the MIP formulation and Ipopt for the rest. 
We ran each instance until convergence or 1 hour had passed, and we recorded the best solution found.

\begin{table}
\centering
\begin{tabular}{r|rrr|rrr}
\toprule
    & \multicolumn{3}{c|}{Shallow} & \multicolumn{3}{c}{Deep} \\
Instance &  MIP &  MPCC &  Embedded &  MIP &  MPCC &  Embedded \\ \midrule
1 &                17.5 &                 \textbf{1.8} &               10.3 &             $3600^*$ &              \textbf{1.6} &            16.7 \\
2 &                22.0 &                 \textbf{1.2} &                8.2 &             $3600^*$ &              \textbf{1.4} &            19.4 \\
3 &                23.9 &                 \textbf{1.5} &                7.5 &             $3600^*$ &              \textbf{1.7} &            14.2 \\
4 &                40.3 &                 \textbf{1.5} &                6.2 &             $3600^*$ &              \textbf{1.4} &             7.2 \\
5 &                19.0 &                 \textbf{1.4} &                4.0 &             $3600^*$ &              \textbf{1.6} &            13.0 \\
6 &                21.0 &                 \textbf{2.1} &                3.0 &             $3600^*$ &              \textbf{1.7} &             7.1 \\
7 &               169.4 &                 \textbf{1.7} &                5.4 &             $3600^*$ &              \textbf{1.2} &            14.1 \\
8 &                14.9 &                 \textbf{1.5} &                5.7 &             $3600^*$ &              \textbf{1.4} &             7.3\\
9 &                86.3 &                 \textbf{1.9} &                4.3 &             $3600^*$ &              \textbf{1.6} &            12.1\\
10 &                15.5 &                 \textbf{1.7} &                5.8 &             $3600^*$ &              \textbf{1.5} &             5.8\\
\bottomrule
\end{tabular}
\caption{Time until a solver finds the best feasible solution for each formulation and for both shallow and deep neural networks on each particular fixing of the oil well problem. The * indicates that the solver could not prove the found solution optimal.}
\label{tab:partial-oil-problem-time}
\end{table}

The results for the experiments performed with the binary variables fixed to arbitrary configurations are presented in Table \ref{tab:partial-oil-problem-time}, which details the solve times, and in Table \ref{tab:partial-oil-problem-obj}, which details the objective values produced.
We observe that, overall, the MPCC formulation clearly outperforms all the other formulations in terms of speed in that it is able to find solutions in about 2 seconds where the embedded network formulations take 10 to 20 seconds and the MIP formulation can take much longer. 
In fact, the MIP formulation fails to find a single solution for any of the instances using the deep neural network configuration, which suggests that optimization problems with deeper networks are more difficult than networks with shallow networks even if the number of nodes is the same. 

\begin{table}
\centering
\begin{tabular}{r|rrr|rrr}
\toprule
    & \multicolumn{3}{c|}{Shallow} & \multicolumn{3}{c}{Deep} \\
Instance &  MIP &  MPCC &  ReLU NN &  MIP &  MPCC &  Embedded  \\ \midrule
1 &              \textbf{1.231} &              \textbf{1.231} &             \textbf{1.231} &              NaN &           \textbf{1.248} &          1.183 \\
2 &              \textbf{1.228} &              \textbf{1.228} &             \textbf{1.228} &            NaN &           \textbf{1.228} &          1.155 \\
3 &              \textbf{1.219} &              \textbf{1.219} &             \textbf{1.219} &            NaN &           1.229 &          \textbf{1.274} \\
4 &              \textbf{1.264} &              \textbf{1.264} &             \textbf{1.264} &            NaN &           \textbf{1.274} &          1.255 \\
5 &              \textbf{1.182} &              \textbf{1.182} &             \textbf{1.182} &            NaN &           1.192 &          \textbf{1.253} \\
6 &              \textbf{1.252} &              \textbf{1.252} &             \textbf{1.252} &            NaN &           1.256 &          \textbf{1.266} \\
7 &              \textbf{1.229} &              \textbf{1.229} &             \textbf{1.229} &            NaN &           \textbf{1.244} &          1.215 \\
8 &              \textbf{1.207} &              \textbf{1.207} &             \textbf{1.207} &            NaN &           1.201 &          \textbf{1.225} \\
9 &              \textbf{1.239} &              \textbf{1.239} &             \textbf{1.239} &            NaN &           1.249 &          \textbf{1.282}\\
10 &              \textbf{1.263} &              \textbf{1.263} &             \textbf{1.263} &            NaN &           \textbf{1.276} &          1.243\\
\bottomrule
\end{tabular}
\caption{Objective value of the best feasible solution for each formulation for both shallow and deep neural networks on each particular fixing of the oil well problem.  NaN indicates that no feasible solution was found.}
\label{tab:partial-oil-problem-obj}
\end{table}

In terms of objective value, for the shallow network all formulations were able to find the global optimum in every instance relatively quickly. 
For the deep networks the different formulations outperformed each other on different instances, leaving no clear winner. 
The MPCC formulation solves each instance the quickest, because  of the relatively small size of each of the neural networks.

Unlike the other two problems, we observed for these instances that the embedded ReLU network formulation did not have trouble with convergence and that, in spite of nondifferentiability, the dual infeasibility was brought down to zero, signifying convergence in all instances. 
We postulate that the reason for the success here as compared with the other instances may be due to the simplicity of the networks: 8 of the networks are single input and single output, and so the difficulties that might emerge in multiple dimensions do not appear.

\section{Conclusions and Future Research}

We have presented three alternative formulations of ReLU deep-neural network constraints as a mixed-integer
problem, an optimization problem with complementarity constraints, and  a problem with the neural network directly embedded. The MIP and MPCC formulations can
be viewed as lifted formulations, and we have shown that the lifting convexifies optimization problems with
deep neural network constraints in the case of the mixed-integer formulation. We have also presented a
warmstart technique that uses training data of the neural network to construct good initial solutions.
We have compared the three formulations on three examples arising in the design of engines, the design
of images that ``fool'' a given classifier, and the assignment of flow in an oil well network. Each formulation has its advantages. We have shown the MIP formulation to be
useful in finding an optimal solution, but  at the cost of a particularly long solve time. 
We observed that the new complementarity constraint formulation generally
outperforms the mixed-integer formulation in terms of solution time but may not find the optimal solution, although it often comes close. 
We also observed that the embedded neural network formulation has the advantage of being scalable and quick to solve but has difficulties with convergence related to the nondifferentiability of the ReLU activation function (which could be rectified by using a smooth activation function like the swish function). 
The experiments reported in this report portray the versatility of each of the different formulations in handling the breadth of optimization problems that may involve neural networks.



\bibliographystyle{abbrv} 
\bibliography{minlp,ml,MPEC}

\vfill
\begin{flushright}
\scriptsize
\framebox{\parbox{\textwidth}{
The submitted manuscript has been created by UChicago Argonne, LLC, Operator of Argonne National Laboratory (“Argonne”). 
Argonne, a U.S. Department of Energy Office of Science laboratory, is operated under Contract No. DE-AC02-06CH11357. 
The U.S. Government retains for itself, and others acting on its behalf, a paid-up nonexclusive, irrevocable worldwide 
license in said article to reproduce, prepare derivative works, distribute copies to the public, and perform publicly 
and display publicly, by or on behalf of the Government.  The Department of Energy will provide public access to these 
results of federally sponsored research in accordance with the DOE Public Access Plan. 
\url{http://energy.gov/downloads/doe-public-access-plan}.
}}
\normalsize
\end{flushright}

\end{document}